\documentclass[a4paper,11pt]{article}
\usepackage[hidelinks]{hyperref}
\hypersetup{
    colorlinks=true,
    linktoc=all,
    linkcolor=red,     %choose some color if you want links to stand out
   citecolor=blue,     %choose some color if you want links to stand out}
}
\usepackage{tikz}
\usepackage{subfigure}
\usepackage{graphicx}
\usepackage{amssymb}
\usepackage{latexsym, bm}
\usepackage{multicol}
\usepackage{indentfirst}
\usepackage{amssymb,amsfonts}
\usepackage{amsmath}
\usepackage{hyperref}
\usepackage{diagbox}
\usepackage[T1]{fontenc}
\usepackage[utf8]{inputenc}
\usepackage{authblk}
\usepackage{setspace}
\renewcommand{\arraystretch}{0.9}

\newtheorem{theorem}{Theorem}[section]
\newtheorem{lemma}[theorem]{Lemma}
\newtheorem{claim}[theorem]{Claim}

\newtheorem{proposition}[theorem]{Proposition}
\newtheorem{definition}[theorem]{Definition}

\newtheorem{conjecture}[theorem]{Conjecture}

\newcommand{\ignore}[1]{}

\begin{document}

\begin{spacing}{1}
\date{}
\title{A step towards the Ramsey-Tur\'{a}n conjecture for $K_3$ and $K_6$}

\author{Xinyu Hu,\footnote{Data Science Institute, Shandong University, Jinan, 250100, P.~R.~ China. Email: {\tt huxinyu@sdu.edu.cn}. }
			\;\; \; Qizhong Lin\footnote{Center for Discrete Mathematics, Fuzhou University, Fuzhou, 350108, P.~R.~China. Email: {\tt linqizhong@fzu.edu.cn}. This research is supported by National Key R\&D Program of China (Grant No. 2023YFA1010202) and NSFC (No.\ 12571361).}
}

\maketitle

\begin{abstract}
Ramsey-Tur\'{a}n  type problems were initiated by Erd\H{o}s and S\'{o}s in 1969.
Given integers $p, q\ge2$, a graph $G$ is $(K_p,K_q)$-free if there exists a red/blue edge coloring of $G$ such that it contains neither a red $K_p$ nor a blue $K_q$. For any $\delta>0$, the Ramsey-Tur\'{a}n number $RT( {n,p,q,\delta n)} $ is the maximum number of edges in an $n$-vertex $(K_p,K_q)$-free graph with independence number at most $\delta n$.
Let $\rho (p, q,\delta ) = \mathop {\lim }\limits_{n \to \infty } \frac{RT(n,p, q,\delta n)}{n^2}$.
Kim, Kim and Liu (2019) showed that $\rho(3,6,\delta)\ge \frac{5}{12}+\frac{\delta}{2}+2\delta^2$ via a skillful construction and conjectured the equality holds for sufficiently small $\delta>0$.
Using Szemer\'{e}di's regularity lemma and a stability argument, we make the first step towards the conjecture by showing that $\rho(3,6,\delta)$ is at most $\frac{5}{{12}} + \frac{\delta }{2}+ 2.1025\delta ^2$.


\end{abstract}

\section{Introduction}
%In this paper, all graphs  we consider are simple, i.e., without loops and multiple edges. Given a graph $G$, $e(G)$ will denote the number of its edges, $v(G)$ the number of its vertices, $\chi(G)$ its chromatic number, $\alpha(G)$ the maximum size of an independent set in it.

Ramsey theorem \cite{ram} implies that for any integers $n_1,\dots,n_k$, there exists a minimum integer, now called the Ramsey number $r=r(n_1,\dots,n_k)$, such that any $k$-coloring of edges of the complete graph $K_r$ contains a $K_{n_i}$ in the $i$th color for some $1\le i\le k$.
Subsequently, Tur\'{a}n posed the problem to determine the maximum number of edges of a $K_{p+1}$-free graph. In particular, Tur\'{a}n \cite{B9,turan54} proved that  the balanced complete $p$-partite graph, known as the {\em Tur\'{a}n graph} $T_{n,p}$(or $T_p(n)$), is the unique extremal graph which attains the maximum number of edges among all $n$-vertex $K_{p+1}$-free graphs.
Since these Tur\'{a}n graphs have large independent sets of size linear in $n$, it is natural to ask for the maximum number of edges of an $n$-vertex $K_{p+1}$-free graph without large independent set. Erd\H{o}s and S\'{o}s \cite{B6} initiated the study of such Ramsey-Tur\'{a}n type problems.

Given integers $p_1,\dots,p_k$, we say that a graph $G$ is $(K_{p_1},\dots,K_{p_k})$-free if there exists a $k$-edge coloring of $G$ with no monochromatic copy of $K_{p_i}$ in the $i$th color for each $1\le i\le k$. The {\em Ramsey-Tur\'{a}n number} $RT(n,p_1,\dots,p_k,m)$ is defined as the maximum number of edges of an $n$-vertex $(K_{p_1},\dots,K_{p_k})$-free graph $G$ with independence number $\alpha(G)\le m$.
Clearly, there is no graph $G$ of order $n$ which is  $(K_{p_1},\dots,K_{p_k})$-free and $\alpha(G)\le m$ if $n\ge r(p_1,\dots,p_k,m)$ from the Ramsey theorem \cite{ram}.
%Indeed, we are interested mainly in the case when $m=o(n)$ but $m/n\to 0$ very slowly, see e.g. the survey by  Simonovits and S\'{o}s \cite{ss}.

\begin{definition}[Ramsey-Tur\'{a}n density]\label{def1}
Given integers $p_1,\dots,p_k$ and $0<\delta<1$, let
\[
\rho (p_1,\dots,p_k,\delta ) = \mathop {\lim }\limits_{n \to \infty } \frac{RT(n,p_1,\dots,p_k,\delta n)}{n^2},
\]
 and $\rho (p_1,\dots,p_k): = \mathop {\lim }\limits_{\delta  \to 0} \rho (p_1,\dots,p_k,\delta )$ is the Ramsey-Tur\'{a}n density of graphs $K_{p_1},\dots,K_{p_k}$.

\end{definition}

%In this paper, we are mainly concerned with the cases where $H_i$ are complete graphs  for $1\le i\le k$. For convenience, we will write $RT(n,{p_1}, \ldots ,{p_k},m)$ instead of $RT(n,K_{p_1}, \ldots ,K_{p_k},m)$, and $\rho({p_1}, \ldots ,{p_k})$ and $\rho({p_1}, \ldots ,{p_k},\delta)$ will be defined similarly.

For $k=1$, the Ramsey-Tur\'{a}n densities of cliques are well understood.
For odd cliques, Erd\H{o}s and S\'{o}s \cite{B6} proved that $\rho( {{{2p + 1}}}) = \frac{1}{2}(1 - \frac{1}{p})$ for all $p \ge 1$. The problem for even cliques is much harder apart from the trivial case $K_2$. Erd\H{o}s and S\'{o}s \cite{B6} showed that $\rho ({{4}}) \le \frac{1}{6}$. As an early application of the regularity lemma,  Szemer\'{e}di \cite{B8} showed that $\rho ({{4}}) \le \frac{1}{8}$. No lower bound on $\rho( {{4}} )$ was known until Bollob\'{a}s and Erd\H{o}s \cite{B1} provided a matching lower bound using an ingenious geometric construction, showing that $\rho ({{4}}) = \frac{1}{8}$. Finally, Erd\H{o}s, Hajnal, S\'{o}s and Szemer\'{e}di \cite{B4} proved $\rho ({{{2p}}}) = \frac{1}{2}( 1 - \frac{3}{3p - 2})$ for all $p \ge 2$.

The quantity $\rho(p,\delta )$ captures more subtle behaviors of the Ramsey-Tur\'an number. Answering two questions posed by Bollob\'{a}s and Erd\H{o}s \cite{B1}, Fox, Loh and Zhao \cite{flz} showed that $\rho(4,\delta )=\frac{1}{8}+\Theta(\delta)$. Recently,  L\"uders and Reiher \cite{lr} determined the exact value of $\rho({p},\delta )$ for each $p\ge3$, in particular,
\[
\left\{ \begin{array}{ll}
\rho({2s+1},\delta )=\frac{1}{2}(\frac{s-1}{s}+\delta) & \textrm{for $s\ge1$,}\vspace{2mm}\\
\rho({2s},\delta )=\frac{1}{2}(\frac{3s-5}{3s-2}+\delta-\delta^2) & \textrm{for $s\geq 2$.}
\end{array} \right.
\]
For more results, we refer the reader to the survey of Simonovits and S\'os \cite{ss} and  some related references \cite{bhs,bls,fs,ll,sud}, etc.

In general, it is much more difficult to determine the exact values of $\rho({p_1}, \ldots,{p_k})$ for $k\ge2$ since the stability structure of a graph that is $(K_{p_1}, \ldots,K_{p_k})$-free is not well understood, which can be reflected from the difficult to determine the Ramsey number $r({p_1}, \ldots,{p_k})$, one can see \cite{li-l,rad} for known small Ramsey numbers. %In particular, for $k\ge3$, we only know $R(3,3,3)=18$ by Greenwood and Gleason \cite{gg} and $R(4,3,3)=30$ by .

 Erd\H{o}s, Hajnal, Simonovits, S\'{o}s and Szemer\'{e}di \cite{B3}  proved that the multicolor Ramsey-Tur\'an density for cliques can be determined by certain weighted Ramsey numbers.
 In particular, they proposed to determine $\rho(p,q)$ for $p,q\ge3$.
 We only know that $\rho(3,q)$ for $q=3,4,5$ and $\rho(4,4)$ from \cite{B3}, and  $\rho(3,6)$ and $\rho(3,7)$ from \cite{H-L}. One can see these values in the following table.

\begin{table}[htb]
\renewcommand{\arraystretch}{1.2}
\begin{center}
\setlength{\tabcolsep}{6mm}{
\begin{tabular}{|c|c|c|c|c|c|}
\hline     \diagbox{$p$}{$\rho(p,q)$}{$q$}    & $3$  & $4$  & $5$ & $6$  & $7$  \\
\hline   $3$ & $\frac{1}{4}$&$\frac{1}{3}$&$\frac{2}{5}$ & $\frac{5}{12}$&$\frac{7}{16}$  \\
\hline   $4$ & $\frac{1}{3}$ &  $\frac{11}{28}$& & &\\
\hline
\end{tabular}}
\end{center}
\caption{The two-colored Ramsey-Tur\'an densities for cliques}
\end{table}

Capturing more subtle behaviors of multicolor Ramsey-Tur\'{a}n numbers, Erd\H{o}s and S\'{o}s \cite{es2} proved in $1979$ that $\rho(3,3,\delta)=\frac{1}{4}+\Theta(\delta)$.
Recently,  Kim, Kim and Liu \cite{kkl} determined that $\rho(3,3,\delta)=\frac{1}{4}+\frac{1}{2}\delta$ for sufficiently small $\delta$, which confirms a conjecture of Erd\H{o}s and S\'{o}s \cite{es2}. Furthermore, they obtained for sufficiently small $\delta>0$, $\rho(3,4,\delta)=\frac{1}{3}+\frac{\delta}{2}+\frac{3\delta^2}{2}$ and  $\rho(3,5,\delta)=\frac{2}{5}+\frac{\delta}{2}$.
They also gave a nice construction showing that
$
 \rho(3,6,\delta)\ge \frac{5}{12}+\frac{\delta}{2}+2\delta^2.
$

\begin{figure}\label{fig-kkl}
  \centering
  \includegraphics[width=9cm]{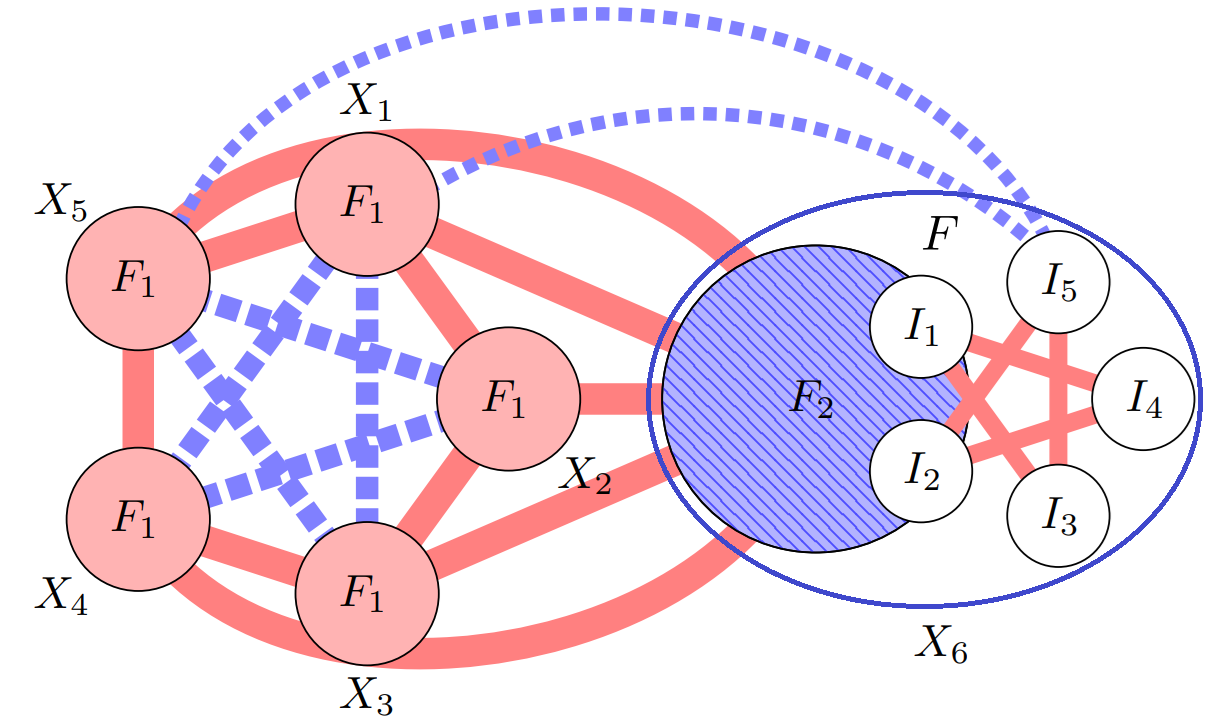}\\\vspace{3mm}
  Fig. \ref{fig-kkl}: Lower bound of $\rho(3,6,\delta)$, where all edges incident to $\cup_{i\in[5]} I_i$ and
$\cup_{i\in[5]} X_i$ are omitted except blue edges between $I_5$ and $X_5\cup X_1$.
  %\caption{}\label{}
\end{figure}

Their construction is as follows (see Fig. \ref{fig-kkl}): Given $d,n$ be integers, let $F(n,d)$ denote an $n$-vertex $d$-regular $K_3$-free graph with independence number $d$. Let $S_\delta  \subseteq ( {0,1} )$ consist of all the rationals $\delta$ for which there exists some $F( {n,d} )$ with $\frac{d}{n} = \delta $. From a result of Brandt \cite{Brandt}, we know that $S_\delta$ is dense in $( {0,\frac{1}{3}} )$.

Assume $6$ divides $n$. In the following, all additions of the subscripts are taken modulo 5. Let $F_1:=F(\frac n6,d_1)$ and $F_2:=F(\frac n6-\frac {3\delta n}2,d_2)$, where $d_i\in[\delta n-o(n),\delta n]$. Let $I=\{v_1,\dots,v_{d_2}\}$ be an independent set of size $d_2$ in $F_2$. Let $I=I_1\sqcup I_2$ be an equipartition of $I$. Let $F$ be an $n/6$-vertex graph obtained from $F_2$ by first adding three clone sets of $I_1$, say $I_i$, with $i\in \{3, 4, 5\}$, and then adding all edges between $I_i$ and $I_{i+2}$ for each $i\in[5]$, and finally adding an additional set of $\frac {3}2(\delta n-d_2)$ isolated vertices. We define $G$ as the graph obtained from $T_6(n)$, by putting a copy of $F$ in $X_6$ and a copy of $F_1$ in $X_i$ for each $i\in[5]$.

Finally, we can see the following
2-edge-coloring of $G$ is $(K_3, K_6)$-free: (1) all edges between $X_i$ and $X_{i+2}$ for each $i\in[5]$ are blue; (2) all edges between $I_i$ and $X_i\cup X_{i+2}$ for each $i\in[5]$ are blue; (3) all edges in $E(G[X_6]\setminus G[\cup_{i\in[5]} I_i])$ are blue; (4) all other edges are red. We know the subgraph induced by all blue edges is $K_3$-free while the subgraph induced by all red edges is $K_6$-free. Now the lower bound follows by noting that $\alpha(G)\le \delta n$ and  $e(G)=(\frac{5}{12}+\frac{\delta}{2}+2\delta^2)n^2$.

\medskip
Based upon the above lower bound, they made the following conjecture.

\begin{conjecture}[Kim, Kim and Liu \cite{kkl}]\label{kkl-2}
For any sufficiently small $\delta>0$,
 $$\rho(3,6,\delta)= \frac{5}{12}+\frac{\delta}{2}+2\delta^2.$$
\end{conjecture}
\medskip

In this paper, we make the first step to Conjecture \ref{kkl-2} by showing the following upper bound.

\begin{theorem}\label{zhu-2}
For any sufficiently small $\delta  > 0$, $ \rho ( {3,6,\delta })\le\frac{5}{{12}} + \frac{\delta }{2}+ 2.1025\delta ^2 $.
\end{theorem}

\noindent
{\bf Notation.} For a graph $G=(V,E)$ with vertex set $V$ and edge set $E$, we use $e(G)$ to denote the number of edges  $|E|$ in $G$. We use $uv$ to denote an edge of $G$. For $X \subseteq V$, we use $e_G(X)$ to denote the number of edges in $X$, and let $G[X]$ denote the subgraph of $G$ induced by $X$. Denote by $\overline X=V\setminus X$ the complement of $X$. For two disjoint subsets $X,Y\subseteq V$, we use $e_G(X,Y)$ to denote the number of edges between $X$ and $Y$.
For a vertex $v\in V$, denote by $N_G(v,X)$ the neighborhood of $v$ in $X$, and $\deg_G(v,X)=|N_G(v,X)|$. In particular,
the neighborhood of a vertex $v$ in $G$ is denoted by $N_G(v)$ and the degree of $v$ in $G$ is $\deg_G(v)=|N_G(v)|$.
$A \sqcup B$ denotes the disjoint union of $A$ and $B$. A complete $p$-partite graph with vertex set $\sqcup_{i=1}^pV_i$, where $|V_i|=n_i$, is denoted by $K_{n_1, \dots, n_p}$. Let $[p]=\{1,2,\ldots,p\}$ and $[p,q]=\{p,p+1,\dots,q\}$, and let $X\choose i$ denote the set of all subsets of a set $X$ of size $i$. We write $a=b\pm c$ if $b-c\le a\le b+c$. We always delete the subscriptions if there is no confusion from the context.
%For undefined terminology, the reader is referred to \cite{bm}.

\medskip\noindent
{\bf Organization:} In Section \ref{pre}, we give useful lemmas for our proofs. In Section \ref{conj-liu}, we will present the proof of Theorem \ref{zhu-2}. Finally, we will mention some interesting problems in Section \ref{crp}. The computation code by LINGO will be placed in the Appendix.

\section{Preliminaries}\label{pre}

Let $G$ be a graph, and let $X,Y\subseteq V(G)$ be disjoint nonempty sets of vertices in a graph $G$. The density of $(X,Y)$ is $d(X,Y)=\frac{e(X,Y)}{|X||Y|}$. For $\varepsilon>0$, the pair $(X,Y)$ is $\varepsilon$-regular in $G$ if for every pair of subsets $X'\subseteq X$ and $Y'\subseteq Y$ with $|X'|\geq \varepsilon |X|$ and $|Y'|\geq \varepsilon |Y|$ we have $| {{d}( {X,Y} ) - {d}( {X',Y'} )} | \le \varepsilon $. Additionally, if $d(X,Y)\geq\gamma$, for some $\gamma >0$, we say that $(X,Y)$ is $(\varepsilon,\gamma)$-regular. Given a $k$-coloring of $E(G)$, we always denote $G_i$ by the spanning subgraph of $G$ induced by all edges of color $i$, where $i\in[k]$.
We say a partition $V(G)=\sqcup_{i=1}^m V_i$ of $G$ is equitable if $||V_i|-|V_j||\le 1$ for all distinct $i$ and $j$. An equitable partition $V(G)= \sqcup_{i=1}^m{{V_i} } $
is said to be an $\varepsilon$-regular partition of a $k$-edge-colored graph $G$ if for each $i\in[m]$, all but at most $\varepsilon m$ choices of $j\in[m]$ satisfy that the pair $(V_i,V_j)$ is $\varepsilon$-regular in $G_\ell$ for each color $\ell\in[k]$.

We will use the following regularity lemma. %For many applications, we refer the reader to Koml\'{o}s and Simonovits \cite{ko-sim} and other related references.
\begin{lemma}[Szemer\'{e}di \cite{sze78}]\label{reg-le}
Suppose $0<\frac{1}{M'}\ll \varepsilon$, $\frac{1}{M}\ll\frac{1}{k}\leq1$, and $n\geq M$. Suppose that $\varphi$ is a $2$-edge-colored of an $n$-vertex graph $G$ and $U_1\sqcup U_2$ is a partition of $V(G)$. Then there exists an $\varepsilon$-regular equitable partition $V(G)= \sqcup_{i=1}^m{{V_i} } $ with $M\leq m\leq M'$ such that for each $i \in[m]$, we have either $V_i \subseteq U_1$ or $V_i \subseteq U_2$.
\end{lemma}

%A useful notion associated with a regular partition is that of a \textbf{reduced graph}.
Given $\varepsilon, \gamma >0$, a graph $G$, a coloring $\varphi : E(G)\rightarrow [2]$, and a partition $V(G)= \sqcup_{i=1}^m{{V_i} } $, we define the {\em reduced graph $R$} as follows: Its vertex set is $V(R)=[m]$, and $ij \in E(R)$ if the pair $(V_i,V_j)$ is $\varepsilon$-regular with respect to $G_k$ for every $k \in [2]$ and $d_{G_k}(V_i,V_j)\geq \gamma$ for some $k \in [2]$.

Given a graph $R$ and $s\ge1$, let $R(s)$ be the {\em blow-up} graph of $R$ obtained by replacing every vertex of $R$ with an independent set of size $s$ and replacing every edge of $R$ with $K_{s,s}$. The following lemma provides some useful properties related to the regular partitions.

\begin{lemma}[Koml\'{o}s and Simonovits \cite{ko-sim}]\label{Count}
 Let $0<\frac{1}{n}\ll \varepsilon \ll \gamma, \frac{1}{h} \leq 1$. Suppose that $H$ is an $h$-vertex graph and $R$ is a graph such that $H\subseteq R(h)$. If $G$ is a graph obtained by replacing every vertex of $R$ with an independent set of size $n$ and replacing every edge of $R$ with an $(\varepsilon , \gamma)$-regular pair, then $G$ contains at least $(\frac{\gamma}{2})^{e(H)}n^{|V(H)|} $ copies of $H$.
\end{lemma}

Applying the degree majorization algorithm used by Erd\H{o}s \cite{e70}, F\"{u}redi \cite{furedi} obtained the following stability result for ${K_{p+1}}$-free graphs.

\begin{lemma}[F\"{u}redi \cite{furedi}]\label{stability}
Let $t$ be a positive integer, and let $G$ be an $n$-vertex ${K_{p+1}}$-free graph with $e( G ) \ge e (T_{n,p}) - t$. Then there exist ${n_1}, \ldots,{n_p}$ such that $e(K_{{n_1}, \ldots ,{n_p}})\geq e (T_{n,p})-2t$ and $|G{\vartriangle}K_{{n_1}, \ldots ,{n_p}}| \leq 3t$. Consequently, $n_i=\frac{n}{p}\pm2\sqrt{t}$ for all $i\in[p]$ and $|G{\vartriangle}T_{n,p}|=O(\sqrt{t}n)$.
\end{lemma}

We will also use the following lemma by Balogh,  Liu and Sharifzadeh \cite[Lemma 3.1]{BLS}, which refines a result of Erd\H{o}s, Hajnal, Simonovits, S\'{o}s and Szemer\'{e}di \cite[Lemma 2]{B3}.

\begin{lemma}[Balogh,  Liu and Sharifzadeh \cite{BLS}]\label{b-indent}
Let $G$ be an $n$-vertex graph with $\alpha(G)\leq {c} n$ for some $0<{c}<1$, and let $\varphi:E(G)\rightarrow [2]$. Then there exists a partition $V(G)=V_1\sqcup V_2$ such that for every $k\in [2]$, $\alpha(G_k[V_k])\leq \sqrt{{c}}n.$
\end{lemma}

 A red-blue edge coloring of $K_5$ is said to be a {\em pentagonlike coloring} if each of the subgraphs induced by all red or blue edges is $C_5$.
The following lemma can be verified directly.
\begin{lemma}\label{jian}
Let $\phi $ be a red-blue edge coloring of $K_5$ such that it contains neither a red $K_3$ nor a blue $K_3$, then $\phi $ must be a pentagonlike coloring.
\end{lemma}

%\begin{lemma}\label{M-s-d}
%Let $(X,Y)$ be an $\varepsilon$-regular pair of density $d$ in a graph $G$. If $B\subseteq Y$ with $|B|\ge \varepsilon|Y|$, then there exist a subset $A\subseteq X$ with $|A|\ge (1-\varepsilon)|X|$ such that each vertices in $A$ is adjacent to at least $(d-\varepsilon)|B|$ vertices in $B$.
%\end{lemma}
%\noindent{\em Proof.} Let $X'$ be the set of vertices with fewer than $(d-\varepsilon)|B|$ neighbors in $B$. Thus $e(X',B)<(d-\varepsilon)|X'||B|$, and so $d(X',B)<d-\varepsilon$. Since $(X,Y)$ is $\varepsilon$-regular, we have that $|X'|<\varepsilon|X|$.

%\hfill$\Box$

%\begin{lemma}[Slicing Lemma \cite{B8}]\label{Sli-L}
%Let $\varepsilon<\alpha,\gamma,\frac{1}{2}$. Suppose that $(X,Y)$ is an $(\varepsilon,\gamma)$-regular pair in a graph $G$. If $X'\subseteq X$ and $Y'\subseteq Y$ satisfies $|X'|\ge \alpha|X|$ and $|Y'|\ge \alpha|Y|$, then $(X',Y')$ is an $(\varepsilon',\gamma-\varepsilon)$-regular pair in $G$, where $\varepsilon':=\max\{\varepsilon/\alpha,2\varepsilon\}$.
%\end{lemma}

\section{Proof of Theorem \ref{zhu-2}}\label{conj-liu}

%We also have the following notation.
Let $G=(V,E)$ be a graph with a partition $ \sqcup_{i=1}^p{{V_i} } $ of $V$. Denote by $G[V_1,\ldots,V_p]$ the $p$-partite subgraph of $G$ induced by $p$-parts $V_1,\ldots,V_p$ of $V$. We say that a partition $ \sqcup_{i=1}^p{{V_i} } $ of $V$ is a {\em max-cut} $p$-partition of $G$ if $e(G[V_1,\ldots,V_p])$ is maximized among all $p$-partition $\sqcup_{i=1}^p{{V_i} } $ of $V$. Denote by $$ \delta^{cr}(G[\sqcup_{i=1}^p{{V_i} } ])= \min_{i,j\in [p],\,v\in V_i}\deg(v,V_j)$$ the minimum crossing degree of $G$ with respect to the partition $\sqcup_{i=1}^p{{V_i} } $.

Given $\varphi:E(G)\rightarrow[2]$, for each $i\in[2]$, we always denote $G_i$ by the spanning subgraph of $G$ induced by all edges of color $i$. We say that $\varphi$ (also $G$) is $(K_{s_1},K_{s_2})$-free if $G_i$ is $K_{s_i}$-free for each $i\in[2]$. We will write $\varphi(A,B)=i$ if $\varphi(e)=i$ for any edge $e$ of $G[A,B]$, and write $\varphi(v,B)$ instead of $\varphi(\{v\},B)$. If $\varphi$ is also defined on $V(G)$, we write $\varphi(A)=i$ if $\varphi(v)=i$ for all $v\in A$. In the following, when we say that a vertex subset $B\subseteq V(G)$ is an independent set, we mean that $B$ spans no edge in $G$. Moreover, all summations of the subscripts are taken modular 5.

\medskip\noindent
{\em Proof sketch of Theorem \ref{zhu-2}.} The weak stability has been proven in our recent paper \cite[Theorem 1.7]{H-L}, stating that the structure of an $n$-vertex $(K_3,K_6)$-free graph $G$ with $\alpha(G)\leq\delta n$ and $e( {{G}} ) = RT( {n,3,6,o( n )} ) + o( {{n^2}} )$ is close to the Tur\'{a}n graph $T_{n,6}$. In order to complete the proof of Theorem \ref{zhu-2}, we need to establish a much more stronger colored stability lemma (Lemma \ref{color stability}).

Let $G$ be an $n$-vertex $(K_3,K_6)$-free graph with $\alpha(G)\leq\delta n$ and $\delta(G) \geq \frac{5n}{6}$. Let
$0<\frac{1}{n}\ll\delta \ll \delta^{*}\ll \frac{1}{m} \ll \varepsilon \ll \gamma\ll 1.$ We first apply Lemma \ref{b-indent} to obtain a partition $V_1^* \sqcup V_2^* $ such that $\alpha(G_k[V_k^*])\leq\delta^{\frac{1}{2}}n$ for $k\in[2]$. Then we apply the regularity lemma (Lemma \ref{reg-le}) with $G, V_1^*, V_2^*, \varepsilon, \varepsilon^{-1}$, and $M'$ to obtain an $\varepsilon$-regular equitable partition $ \sqcup_{i=1}^m{{V_i} } $ with $\varepsilon^{-1} \leq m \leq M'$ which refines the partition $V_1^* \sqcup V_2^* $.

We will show that for any 2-edge coloring of $G$, there exists a partition of $V(G)=\sqcup_{i=1}^6X_i$ that satisfies $8$ properties (see Lemma \ref{color stability}). To this end, we first use the structure of the reduced graph for the regular pairs to obtain the edges between parts $X_i$ and $X_{i+1}$ for $i\in[5]$ are almost in color 1 (see Proposition \ref{J'}), then we use these edges distributions to obtain that $\deg_1(x,X_j)$ are small for any $k\in[5]$, and $x\in X_k$, and $j\in [k-1,k+1]$ (see Claim \ref{v1-4} and Claim \ref{v1-3}). Moreover, we can also obtain  $\deg_2(x,X_j)$ are small for any $k\in[5]$, and $x\in X_k$, and $j\in \{k-2,k,k+2\}$ (see Claim \ref{v1-2} and Claim \ref{y-X1}). Finally, we can verify these 8 properties of the colored stability lemma.

Finally, on the contrary, we may assume that $e(G)$ is greater than the upper bound from Theorem \ref{zhu-2}. Then, we apply the lemma \cite[Lemma \ref{kkl-c}]{kkl} to obtain a subgraph $G'$ with a lower bound of $e(G')$. Subsequently, we apply the colored stability lemma for $G'$ to obtain $V(G')=\sqcup_{i=1}^6X_i$
that satisfies 8 properties. In particular, in $X_6$, we define five small disjoint subsets $I_1, \ldots, I_5$, and let $I=\sqcup_{i=1}^5I_i$.
Finally, we separate the proof into two cases according to that whether all edges in $G'[I,X_6\setminus I]$ are colored 1. For each case, we can get a contradiction.

\medskip
The details of the proof are as follows.

\subsection{Stability of the colored graph}

\begin{lemma}\label{color stability}
Suppose $0<\frac{1}{n} \ll \delta \ll \gamma \ll 1$. Let $G$ be an $n$-vertex $(K_3,K_6)$-free graph with $\alpha(G)\leq\delta n$ and $\delta(G) \geq \frac{5n}{6}$. Then, for any $2$-edge coloring $\varphi: E(G)\rightarrow [2]$, there exists a partition $ \sqcup_{i=1}^6{{X_i} } $ of $V(G)$ such that the following hold (relabel these $X_i$'s if necessary).
\begin{description}
  \item[($P_1$)] For each $i\in [6]$, we have $|X_i|=\frac{n}{6}\pm 2\gamma^{\frac{1}{4}}n.$
  \item[($P_2$)] There exists some part, say $X_6$, such that $\alpha(G_1[X_6])\leq \gamma^{\frac{1}{4}}n.$
  \item[($P_3$)] For each vertex $v\in X_6$, $\min \{{\deg_{G_1}(v,X_i), \deg_{G_1}(v,X_{i+2})\} \leq \gamma^\frac{1}{59}}n$ for $i\in [5]$.
  \item[($P_4$)] For each vertex $v\in X_6$, $\min\limits_{i\in [5]} \{{\deg_{G_1}(v,X_i\sqcup X_{i+1})\} \leq \gamma^\frac{1}{60}}n$.
  \item[($P_5$)] For each $i\in[6]$, we have $\Delta(G[X_i])\leq \gamma^{\frac{1}{117}}n.$
  \item[($P_6$)] $\delta^{cr}(G[ \sqcup_{i=1}^6{{X_i} } ])\geq \frac{n}{6}-\gamma^{\frac{1}{118}}n$.
  \item[($P_7$)] For each vertex $v \in X_i$ with $i \in [5]$, $\deg_{G_2}(v,X_6)\geq| {{X_6}} | - {\gamma ^{\frac{1}{{119}}}}n.$
  \item[($P_8$)]  For each vertex $v \in X_i$ with $i \in [5]$, $\alpha(G_2[X_i])\leq \gamma^{\frac{1}{4}}n$, and $\deg_{G_1}(v,X_{j_1})\geq |X_{j_1}|-\gamma^{\frac{1}{119}}n$ where $j_1\in \{i-2,i+2\}$, and $\deg_{G_2}(v,X_{j_2})\geq |X_{j_2}|-\gamma^{\frac{1}{119}}n$ where $j_2\in \{i-1,i+1\}$.

\end{description}
\end{lemma}
{\bf Proof.}
We choose the parameters as follows:
$0<\frac{1}{n}\ll\delta \ll \delta^{*}\ll \frac{1}{m} \ll \varepsilon \ll \gamma\ll 1.$ We first apply Lemma \ref{b-indent} with ${c=\delta}$ to obtain a partition $V_1^* \sqcup V_2^* $ such that $\alpha(G_k[V_k^*])\leq\delta^{\frac{1}{2}}n$ for $k\in[2]$. Then we apply Lemma \ref{reg-le} with $G, V_1^*, V_2^*, \varphi, \varepsilon, \varepsilon^{-1}$, and $M'$ playing the roles of $G, U_1, U_2, \varphi, \varepsilon, M$, and $M'$ to obtain an $\varepsilon$-regular equitable partition $ \sqcup_{i=1}^m{{V_i} } $ with $\varepsilon^{-1} \leq m \leq M'$ which refines the partition $V_1^* \sqcup V_2^* $. For convenience, we assume that  $|V_i|=\frac{n}{m}$ for $1\le i\le m$.
Let $R$ be its reduced graph defined on $[m]$. From \cite[Theorem 1.7]{H-L}, we know $|G\Delta T_{n,6}|\leq \delta^{*}n^{2}$, and so the number of $K_7$ in $G$ is at most $\delta^{*}n^{7}$. Therefore,
\begin{align}\label{str-R}
 \text{$R$ is $K_7$-free,}
\end{align}
since otherwise Lemma \ref{Count} implies that $G$ contains at least $\frac{1}{2}(\frac{\gamma}{2})^{21}(\frac{n}{m})^7>\delta^*n^7$ copies of $K_7$, a contradiction.
Moreover,
\begin{align}\label{str-dt}
 \delta(R) \geq \left(\frac{5}{6}-3\gamma\right)m,
\end{align}
 since otherwise $\delta(G)\le (\delta(R)+\varepsilon m+2\gamma m+1)\frac{n}{m}\le \frac{5}{6}n-(\gamma-\varepsilon-\frac{1}{m})n< \frac{5}{6}n$ due to $\frac{1}{m} \ll \varepsilon \ll \gamma$, a contradiction  again.  Thus, by Lemma \ref{stability},
\begin{equation}\label{4.1}
  |R\Delta T_{m,6}|\leq \gamma^\frac{1}{3}m^2.
\end{equation}

We define a coloring $\phi^R:V(R) \sqcup E(R)\rightarrow [2]$, induced by $\varphi$, as follows:
\begin{description}
\item[(i)] for each $j \in [m]$, we have $\phi^R(j)=k$ if $V_j\subseteq V_k^*$ for $k\in[2]$;
\item[(ii)] for each edge $pq \in E(R)$, we have $\phi^R(pq)=1$ if $d_{G_1}(V_p,V_q)\geq\gamma$, and we have $\phi^R(pq)=2$ if $d_{G_1}(V_p,V_q)<\gamma$ and $d_{G_2}(V_p,V_q)\geq\gamma.$ This is reasonable from the definition of the reduced graph.
\end{description}

For each $pq\in E(R)$, define the weight of $pq$ by $\omega(pq)=d_{G_1}(V_p,V_q)$ \textbf{iff} $d_{G_1}(V_p,V_q)\geq\gamma$, and  $\omega(pq)=d_{G_2}(V_p,V_q)$ \textbf{iff} $d_{G_1}(V_p,V_q)<\gamma$ but $d_{G_2}(V_p,V_q)\geq\gamma$.
%we let $\omega(pq)=d_{G_{\phi^R(pq)}}(V_p,V_q)$ be the weight on $E(R)$,
Now, we consider $R$ as a weighted graph. Note that for each fixed $p\in[m]$, {$|V_p|=|V_q|=\frac{n}{m}$, and for each $v\in V_p$ the neighbors of $v$ in $V_p$ is at most $\frac{n}{m}$}, and all but at most $\varepsilon m$ choices of $q\in[m]$ satisfy that the pair $(V_p,V_q)$ is $\varepsilon$-regular, and  for each non-edge $pq$ the density of $(V_p,V_q)$ in $G_k$ is at most $\gamma$, so we have that
\begin{equation}\label{4.2}
  \sum\limits_{q \in {N_R}( p )} {\omega( {pq} )}  \ge {\left\{{\frac{n}{m}}\cdot\left(\delta ( G ) -{\frac{n}{m}}\right)-
  (\varepsilon m  + m\cdot2\gamma ) \left( {\frac{n}{m}}\right)^2 \right\}\bigg/{\left( {\frac{n}{m}} \right )^2}} \ge \left( {\frac{5}{6} - 3\gamma } \right)m.
\end{equation}

Let $R'$ be the subgraph obtained from $R$ by deleting all edges of weight at most $\frac{1}{2}+\gamma$. We claim that
\begin{equation}\label{4.4}
  e(R)-e(R')\leq \gamma^\frac{1}{4}m^2.
\end{equation}
Indeed, for each $p \in V(R)$, we have
\[
\sum\limits_{q \in {N_R}( p )} {\omega( {pq} )}  \le 1\cdot {\deg_{R'}}( p ) + \left( {\frac{1}{2} + \gamma } \right)( {{\deg_R}( p ) - {\deg_{R'}}( p )} )\le\frac{1}{2}{\deg_{R'}}( p )+\frac{1}{2}{\deg_{R}}( p )+\gamma m,
\]
which together with (\ref{4.2}) yield that
\begin{equation}\label{4.3}
  \deg_{R'}(p)\geq\frac{5m}{3}-\deg_R(p)-8\gamma m.
\end{equation}

 Since $e(R)\leq \frac{5}{12}m^2+\gamma^\frac{1}{3}m^2$ from (\ref{4.1}), and $\delta(G)\geq \frac{5n}{6}$, we have that

\begin{align*}
  \frac{5}{{12}}{n^2} \le e( G ) &\le e( {R'} )\frac{{{n^2}}}{{{m^2}}} + ( {e( R ) - e( {R'} )} )\left( \frac{1}{2} + \gamma  \right)\frac{{n^2}}{{{m^2}}} +  \left( {2\gamma  + \varepsilon  + \frac{1}{m}} \right){n^2}
  %\\& =(e( {R'} )-e(R))\frac{{{n^2}}}{{{m^2}}} +e(R)\frac{{{n^2}}}{{{m^2}}}+ ( {e( R ) - e( {R'} )} )\left( \frac{1}{2} + \gamma  \right)\frac{{n^2}}{{{m^2}}} +  \left( {2\gamma  + \varepsilon  + \frac{1}{m}} \right){n^2}
  \\
   &=( e( R ) - e( {R'} ))\left(\gamma- \frac{1}{2}\right)\frac{{n^2}}{{m^2}}+e(R)\frac{{{n^2}}}{{{m^2}}}+\left( {2\gamma  + \varepsilon  + \frac{1}{m}} \right){n^2}\\
   &\le %( e( R ) - e( {R'} ))(\gamma- \frac{1}{2})\frac{{n^2}}{{m^2}}+\left(\frac{5}{12}m^2+\gamma^\frac{1}{3}m^2\right)\frac{{n^2}}{{m^2}}+\left( {2\gamma  + \varepsilon  + \frac{1}{m}} \right){n^2}\\
   ( e( R ) - e( {R'} ))\left(\gamma-\frac{1}{2}\right)\frac{{n^2}}{{m^2}}+
   \left(\frac{5}{12}n^2+\gamma^\frac{1}{3}n^2\right)+3\gamma{n^2}
\end{align*}
by noting (\ref{4.1}) and $\frac{1}{m} \ll \varepsilon \ll \gamma$. Therefore,
$( e( R ) - e( {R'} ))(\frac{1}{2}-\gamma)\frac{{n^2}}{{m^2}}\le(\gamma^\frac{1}{3}+ 3\gamma){n^2},$
and so
$e(R)-e(R')\leq \gamma^\frac{1}{4}m^2$ follows as desired.

%We will omit $\gamma$ in the term ``$\gamma$-generalized clique."
Given a weighted graph $R$ with weight $\omega : E(R)\rightarrow(0,1]$ and $Y\subseteq X\subseteq V(R)$, a {\em $\gamma$-generalized clique} $Z_t$ of order $t=|X|+|Y|$ on $(Y,X)$ is a clique on $X$ with $\omega(e)>\frac{1}{2}+\gamma$ for every edge $e$ in $Y$.
For each $k\in[2]$ and $Y\subseteq X\subseteq V(R)$, we say that a {\em $\gamma$-generalized clique} $Z_t$ in $R$ on $(Y,X)$ is on color $k$ if $\phi^R(j)=\phi^R(pq)=k$ for all $j \in Y$ and $pq \in {X\choose 2}$. We say that $R$ is $(Z_{t_1}, Z_{t_2})$-free if there is no $Z_{t_k}$ of color $k$ for any $k\in[2]$.

 Since $G$ is $(K_3,K_6)$-free, by noting \cite[Claim 3.2]{H-L}, we have that
\begin{align}\label{Z3Z6-free}
\text{$R$ is $(Z_3,Z_6)$-free.}
\end{align}

Let $\sqcup_{i=1}^6  U_i$ be a {\bf max-cut $6$-partition} of $R$. The desired partition of $V(G)$ will be adjustment of this partition. By (\ref{4.1}) and Lemma \ref{stability}, we obtain that
\begin{equation}\label{4.5}
 \sum\limits_{i \in [6]} {e( {R[ {{U_i}} ]} )}  \le {\gamma ^{\frac{1}{3}}}{m^2}, ~~~~| {{U_i}} | = \frac{m}{6} \pm {\gamma ^{\frac{1}{4}}}m,
\end{equation}
and from (\ref{str-dt}),
\begin{equation}\label{4.6}
  {\delta ^{cr}}( {R[  \sqcup_{i=1}^6{{U_i} } ]} ) \ge   \frac{1}{2}\left({\delta ( R ) - 4\mathop {\max }\limits_{i \in [ 6 ]} | {{U_i}} |}\right) \ge \frac{m}{13}.
\end{equation}
\medskip

We will discuss the color patterns of vertices and edges of $R$ in $\phi^R$. First, we show that each vertex set $U_i$ is monochromatic for $i\in[6]$ in $\phi^R$.

\begin{claim}\label{mono}
For every $i\in[6]$, there is $k\in[2]$ such that $\phi^R(U_i)=k$ and
$
\alpha ( {{G_k}[ {\sqcup_{\ell \in {U_i}} {{V_\ell}} } ]} ) \le \sqrt \delta  n.
$
\end{claim}
{\bf Proof.}~On contrary, we may assume that $\phi^R(U_1)\neq k$ for any $k \in [2]$ without loss of generality. We first show the following propositions.

\begin{proposition}\label{two-cases}
There exist vertices $u,v \in U_1$, and $v_i \in U_i$ for $i\in [2,6]$ such that $\phi^R(u)=1$ and $\phi^R(v)=2$, and  $\{{v_2,\dots,v_6}\}$ induces a $K_5$ in $R'$. Furthermore, one of the following holds:
\smallskip

(i) $vv_i \in E(R')$, $uv_i \in E(R)$ for each $i\in [2,6]$.

\smallskip

(ii) $uv$, $uv_i$, $vv_i \in E(R)$ for each $i\in [2,6]$.
\end{proposition}
{\bf Proof.}
Let $U_1(1)= \{ {v \in {U_1}:{\phi ^R}( v ) = 1} \}$ and $U_1(2)= \{ {v \in {U_1}:{\phi ^R}( v ) = 2} \}$. Clearly, $U_1(1)$ and $U_1(2)$ form a partition of $U_1$.

First, suppose $|U_1(2)|\geq \frac{m}{100}$. Fix an arbitrary vertex $u\in U_1$ with $\phi^R(u)=1$. Then, more than half of the vertices $v$ in $U_1(2)$ satisfy ${\deg_{R'}}( v ) \ge {\deg_R}( v ) - {\gamma ^{\frac{1}{5}}}m$ since $ e(R)-e(R')\leq \gamma^\frac{1}{4}m^2$ from (\ref{4.4}). Also,  more than half of the vertices $v$ in $U_1(2)$ satisfy $|N_{R}(v,U_1)|\leq \gamma^\frac{1}{4} m$ since $\sum\limits_{i \in [6]} {e( {R[ {{U_i}} ]} )}  \le {\gamma ^{\frac{1}{3}}}{m^2}$ from (\ref{4.5}). Hence there exists some vertex $v\in U_1(2)$ such that for each $i\in [2,6]$,
\[
| {{N_{R'}}( {v,{U_i}} )} | \ge \delta ( R ) - 4\mathop {\max }\limits_{i \in [ 6 ]} | {{U_i}} | - {\gamma ^{\frac{1}{5}}}m- {\gamma ^{\frac{1}{4}}}m \ge \frac{m}{7}.
\]
 This and the fact ${\delta ^{cr}}( {R[ \sqcup_{i=1}^6{{U_i} } ]} ) \ge \frac{m}{13}$ from (\ref{4.6}) yield that for each $i\in [2,6]$,
\[
| {{N_R}( {u,{U_i}} ) \cap {N_{R'}}( {v,{U_i}} )} | \ge \frac{m}{7} + \frac{m}{13} - | {{U_i}} | \ge \frac{m}{{20}}.
\]
Together with (\ref{4.1}) and (\ref{4.4}), we can find some vertex $v_i\in {{N_R}( {u,{U_i}} ) \cap {N_{R'}}( {v,{U_i}} )}$ for $i\in[2,6]$ such that $\{{v_2,\dots,v_6}\}$ induces a $K_5$ in $R'$, yielding Proposition \ref{two-cases} (i).

We now assume $|U_1(2)|< \frac{m}{100}$. Fix an arbitrary vertex $v\in U_1(2)$. If $|N_R(v,U_1)| > \frac{m}{50}$, then $|N_R(v,U_1(1))|\geq\frac{m}{100}$. As $\sum\limits_{i \in [6]} {e( {R[ {{U_i}} ]} )}  \le {\gamma ^{\frac{1}{3}}}{m^2}$, more than half of the vertices $u$ in $N_R(v,U_1(1))$ satisfy $|N_R(u,U_1)|\leq \gamma^\frac{1}{4} m$. Thus there is some vertex $u \in N_R(v,U_1(1))$ with $|N_R(u,U_i)|\geq \frac{m}{7}$ for each $i\in [2,6]$.  This together with (\ref{4.6}) yields that for each $i\in [2,6]$,
\[
| {{N_R}( {u,{U_i}} ) \cap {N_{R}}( {v,{U_i}} )} | \ge \frac{m}{7} + \frac{m}{13} - | {{U_i}} | \ge \frac{m}{{20}}.
\]
Therefore, again by (\ref{4.1}) and (\ref{4.4}), there exists $v_i\in {{N_R}( {u,{U_i}} ) \cap {N_{R}}( {v,{U_i}} )}$ for all $i\in[2,6]$, such that $\{{v_2,\dots,v_6}\}$ induces a $K_5$ in $R'$, yielding Proposition \ref{two-cases} (ii).
Thus we may assume that $|N_R(v,U_1)| \leq \frac{m}{50}$, and so ${\deg_R}( v ) \le \sum\limits_{i = 2}^6 {| {{U_i}} |}  + \frac{m}{{50}} \le ( {\frac{5}{6} + \frac{1}{{40}}} )m.$
This together with  (\ref{4.3}) yield
\[
{\deg_{R'}}( v ) \ge \frac{5}{3}m - \left( {\frac{5}{6} + \frac{1}{{40}}} \right)m - 8\gamma m \ge \left( {\frac{5}{6} - \frac{1}{{30}}} \right)m.
\]
Hence, for each $i \in [2,6]$,
\[
| {{N_{R'}}( {v,{U_i}} )} | \ge {\deg_{R'}}( v ) - | {{N_R}( {v,{U_1}} )} | - 4\mathop {\max }\limits_{i \in [ 6 ]} | {{U_i}} | \ge \left( {\frac{1}{6} - \frac{1}{{30}} - \frac{1}{{40}}} \right)m \ge \frac{m}{{10}}.
\]
By (\ref{4.6}), there exists a vertex $u\in U_1(1)$ such that for each $i \in [2,6]$,
\[
| {{N_R}( {u,{U_i}} ) \cap {N_{R'}}( {v,{U_i}} )} | \ge \frac{m}{{10}} + \frac{m}{13} - | {{U_i}} | \ge \frac{m}{{100}}.
\]
Together with (\ref{4.1}) and (\ref{4.4}), we have that there exists $v_i\in {{N_R}( {u,{U_i}} ) \cap {N_{R'}}( {v,{U_i}} )}$ for $i\in[2,6]$ such that $\{{v_2,\dots,v_6}\}$ induces a $K_5$ in $R'$, yielding Proposition \ref{two-cases} (i).\hfill$\Box$

\medskip

Let $u,v \in U_1$, and $v_i \in U_i$ for $i\in [2,6]$ be chosen as in Proposition \ref{two-cases}.
Since  $\phi^R(u)=1$ from Proposition \ref{two-cases},  we must have $\phi^R(uv_i)=2$ for $i \in [2,6]$; otherwise, we get a $Z_3$ of color $1$, contradicting the fact that $R$ is $(Z_3,Z_6)$-free.

\begin{proposition}\label{color 2}
For $i \in [2,6]$, we have $\phi^R(v_i)=2$.
\end{proposition}
{\bf Proof.} On contrary, suppose that $\phi^R(v_2)=1$ without loss of generality. Then, for each $i\in[3,6]$, we have $\phi^R(v_2v_i)=2$ since $R$ is $(Z_3,Z_6)$-free. Furthermore, for distinct $j',j\in [3,6]$, there exists at least one edge $v_{j'}v_j$ with $\phi^R(v_{j'}v_j)=1$; otherwise, we get a $Z_6$ of color $2$ on $(\emptyset,\{u,v_2,v_3,\dots,v_6\})$, contradicting the fact that $R$ is $(Z_3,Z_6)$-free. By symmetry, we may assume that $\phi^R(v_3v_4)=1$, then $\phi^R(v_3)=\phi^R(v_4)=2$ since $R$ is $(Z_3,Z_6)$-free.

Suppose now that $\phi^R(v_5v_6)=2$, then we must have $\phi^R(v_4v_5)=1$ or $\phi^R(v_4v_6)=1$; otherwise, $(\{v_4\},\{u,v_2,v_4,v_5,v_6\})$ would be a $Z_6$ of color $2$, which in turn implies that $\phi^R(v_5)=2$ or $\phi^R(v_6)=2$ to avoid a $Z_3$ of color $1$. Then, $\phi^R(v_3v_5)=2$ or $\phi^R(v_3v_6)=2$ again to avoid a $Z_3$ of color $1$. However then $(\{v_3,v_5\},\{u,v_2,v_3,v_5\})$ or $(\{v_3,v_6\},\{u,v_2,v_3,v_6\})$ would form a $Z_6$ of color $2$, a contradiction.

Hence, we may assume that $\phi^R(v_5v_6)=1$, and so $\phi^R(v_5)=\phi^R(v_6)=2$ to avoid a $Z_3$ of color $1$, which in turn implies that $\phi^R(v_3v_5)=\phi^R(v_3v_6)=1$; otherwise, $(\{v_3,v_5\},\{u,v_2,v_3,v_5\})$ or $(\{v_3,v_6\},\{u,v_2,v_3,v_6\})$ would form a $Z_6$ of color $2$. However then we shall get a $Z_3$ of color $1$ on $(\emptyset,\{v_3,v_5,v_6\})$, a contradiction.\hfill$\Box$
\medskip

Now we continue the proof of Claim \ref{mono}.
Suppose Proposition \ref{two-cases} (i) holds.
Recall that $\{{v_2,\dots,v_6}\}$ forms a $K_5$ in $R'$. Clearly, this $K_5$ contains no $K_3$ of color $1$. We claim that this $K_5$ contains no $K_3$ of color $2$ too. On the contrary, by symmetry, we may assume that $\{v_2,v_3,v_4\}$ forms a $K_3$ of color $2$. Since $\phi^R(v_i)=2$ for $i\in[2,4]$ from Proposition \ref{color 2}, and $\omega(v_iv_j)\geq \frac{1}{2} +\gamma$ for $2\le i<j\le4$ from the definition of $R'$, we have that $(\{v_2,v_3,v_4\},\{v_2,v_3,v_4\})$ forms a $Z_6$ of color $2$, a contradiction. Therefore,  by Lemma \ref{jian}, the edge coloring of this $K_5$ must be a \emph{pendagonlike coloring}.
Without loss of generality, we may assume that all edges of the cycle $C_5:v_2v_3v_4v_5v_6v_2$ receive color $1$.  Suppose that $\phi ^R(vv_2)=1$, then $\phi^R(vv_3)=\phi^R(vv_6)=2$; otherwise, $\{v,v_2,v_3\}$ or $\{v,v_2,v_6\}$ would form a $Z_3$ of color $1$. However, then we have a $Z_6$ on $(\{v,v_3,v_6\},\{v,v_3,v_6\})$  of color $2$, a contradiction. So we may assume $\phi ^R(vv_2)=2$, then $\phi^R(vv_4)=\phi^R(vv_5)=1$; otherwise, $\{v,v_2,v_4\}$ or $\{v,v_2,v_5\}$ would yield a $Z_6$ of color $2$. However then we have a $Z_3$ of color $1$ on $\{v,v_4,v_5\}$, again a contradiction.

Now suppose Proposition \ref{two-cases} (ii) holds. By Proposition \ref{color 2} and $\{{v_2,\dots,v_6}\}$ induces a $K_5$ in $R'$ and $R$ is $(Z_3,Z_6)$-free, then the edge coloring of this $K_5$ must be {\em pentagonlike coloring}. We may assume all edges of the cycle $C_5:v_2v_3v_4v_5v_6v_2$ receive color $1$. Thus, we must have $\phi ^R(vv_2)=1$ or $\phi ^R(vv_4)=1$; otherwise, $(\{v_2,v_4\},\{u,v,v_2,v_4\})$ would be a $Z_6$ of color $2$, which in turn implies that $\phi^R(vv_6)=2$ or $\phi^R(vv_5)=2$; otherwise, $\{v,v_6,v_2\}$ or $\{v,v_5,v_4\}$ would form a $Z_3$ of color $1$. So, we get $\phi^R(vv_3)=1$; otherwise, $(\{v_3,v_6\},\{v,u,v_3,v_6\})$ or $(\{v_3,v_5\},\{v,u,v_3,v_5\})$ would be a $Z_6$ of color $2$. But then $\{v,v_3,v_2\}$ or $\{v,v_3,v_4\}$ would yield a $Z_3$ of color $1$, a contradiction.

Therefore, for each $i \in [6]$, there exists $k\in[2]$ such that $\phi^R(U_i)=k$. Then by the definition of $\phi^R$, we have $\sqcup_{\ell \in {U_i}} {{V_\ell}}  \subseteq V_k^*$, and so
$\alpha ( {{G_k}[ {\sqcup_{\ell \in {U_i}} {{V_\ell}} } ]} ) \le \alpha ( {{G_k}[ {V_k^*} ]} ) \le \sqrt \delta  n$ as desired. This completes the proof of Claim \ref{mono}.
\hfill$\Box$

\medskip

Recall that $R'$ is the graph obtained from $R$ by deleting all edges of weight at most $\frac{1}{2}+\gamma$ and $\sqcup_{i=1}^6  U_i$ is a max-cut $6$-partition of $R$, then $\sqcup_{i=1}^6  U_i$ is also a $6$-partition of $R'$.
Define
$$R_0=\left\{v\in R:\deg_R(v)-\deg_{R'}(v)\geq\gamma^\frac{1}{12}m\right\}, \; \text{and} \; \;
{R_1} = \bigcup_{i = 1}^6 {\left\{ {v \in {U_i}:{\deg_R}( {v,{U_i}} ) \ge {\gamma ^{\frac{1}{{12}}}}m} \right\}} .$$
Since $ e(R)-e(R')\leq \gamma^\frac{1}{4}m^2$ from (\ref{4.4}), we have
$$ 2{\gamma ^{\frac{1}{4}}}{m^2}  \ge 2( {e( R ) - e( {R'} )} )
    = \sum\limits_{v \in R} {( {{\deg_R}( v ) - {\deg_{R'}}( v )} )}
   \ge \sum\limits_{v \in {R_0}} {( {{\deg_R}( v ) - {\deg_{R'}}( v )} )}
    \ge | {{R_0}} |{\gamma ^{\frac{1}{12}}}m,$$
which implies that $$|R_0|\leq 2\gamma^\frac{1}{6} m.$$
Similarly, since $\sum\limits_{i \in [6]} {e( {R[ {{U_i}} ]} )}  \le {\gamma ^{\frac{1}{3}}}{m^2}$ from (\ref{4.5}), we have $$|R_1|\leq 2\gamma^\frac{1}{4} m.$$
For each $i\in[6]$, define
$$W_i = {U_i}\backslash ({R_0} \sqcup R_1 ).$$
 Note that $|U_i|=\frac{m}{6}\pm \gamma^{\frac14}m$, so we obtain that for each $i\in[6]$,
\begin{equation}\label{4.7}
|W_i|= \frac{m}{6}\pm\gamma^\frac{1}{7}m.
\end{equation}
From the definition of $R_1$, each vertex of $W_i$ has at most $\gamma^{\frac{1}{{12}}}m$ neighbors in its own set. Therefore, the minimum crossing degree
\begin{equation}\label{4.8}
  {\delta ^{cr}}( {R'[  \sqcup_{i=1}^6{{W_i} }  ]} ) \ge (\delta ( R ) - {\gamma ^{\frac{1}{{12}}}}m) - {\gamma ^{\frac{1}{{12}}}}m- (|R_0|+|R_1|)- 4\mathop {\max }\limits_{i \in [ 6 ]} | {{U_i}} | \ge \frac{m}{6} - {3\gamma ^{\frac{1}{{12}}}}m.
\end{equation}

Now we will consider the color patterns of edges between $U_i$ and $U_j$ for $1\le i<j\le 6$.
Let $$R''=R'[ \sqcup_{i=1}^6{{W_i} }],$$ and denote $\phi^{R''}$ by the two-coloring of vertices and edges restricted on $R''$, that is ${\phi ^{R''}} := { {{\phi ^R}} |_{R''}}.$

\begin{claim}\label{1-v-c}
Relabel $U_1,\dots, U_6$ if necessary, we have

 \smallskip
(i) $\phi^R(U_6)=1$, and for each other $i \in [5]$, $\phi^R(U_i)=\phi^{R}(U_i,U_6)=2$.

 \smallskip
(ii) $\phi^{R''}(W_i,W_{i+1})=2$, and $\phi^{R''}(W_i,W_{i+2})=1$  for $i \in [5]$.
\end{claim}
{\bf Proof.}~ (i) Suppose to the contrary that $\phi^R(U_i)=2$ for all $i\in[6]$. Then, it follows by $(\ref{4.7})$ and $(\ref{4.8})$ that for each $i\in [6]$, there exists $v_i\in W_i$ such that $\{v_1,\dots ,v_6\}$ induces a $K_6$ in $R''$. Note that $r(3,3)=6$. If there exists a triangle of color $1$, then we have a $Z_3$ of color $1$, a contradiction. If there exists a triangle of color $2$, then we obtain a $Z_6$ of color $2$ since $\phi^R(U_i)=2$ for all $i\in[6]$ and all edges in $R''$ have weights at least $\frac12+\gamma$, again a contradiction. Thus, we may assume that $\phi^R(U_6)=1$ since all vertices receive the same color by Claim \ref{mono}. We must have $\phi^{R}(U_i,U_6)=2$ for $i\in[5]$ since otherwise $R$ would contain a  $Z_3$ of color $1$. Note that $\phi^R(v_6)=1$ and $\phi^R(v_iv_6)=2$ for $i\in[5]$. By a similar argument as Proposition \ref{color 2}, we can show that $\phi^R(v_i)=2$ for all $i\in[5]$. Thus Claim \ref{mono} implies that $\phi^R(U_i)=2$ for each $i \in [5]$.

\medskip
(ii)~We first show that for distinct $i,j\in[5]$, $\phi^{R''}(W_i,W_{j})=1$ or $\phi^{R''}(W_i,W_{j})=2$.
 On the contrary, suppose that there exist $w_{11}, w_{12} \in W_1$ and $w_{21}, w_{22} \in W_2$ such that $\phi^{R''}(w_{11}w_{21})=1$ but $\phi^{R''}(w_{12}w_{22})=2$ without loss of generality (where $w_{11}=w_{12}$ or $w_{21}=w_{22}$ is possible).
For each $i\in[3,5]$, since ${\delta ^{cr}}( {R'[  \sqcup_{i=1}^6{{W_i} }  ]} )\ge \frac{m}{6} - {3\gamma ^{\frac{1}{{12}}}}m$ from $(\ref{4.8})$,  we have
 \[
| {{N_{R'}}( {\{ {w_{11},w_{12},w_{21},w_{22}} \},{W_i}} )} | \ge 4{\delta^{cr}}( {R'[ \sqcup_{i=1}^6{{W_i} }   ]} ) - 3| {{W_i}} | \ge \frac{m}{6} - {13\gamma ^{\frac{1}{{12}}}}m.
\]
Then, we can choose $v_i\in N_{R'}( \{ {w_{11},w_{12},w_{21},w_{22}})$ for $i\in[3,5]$ such that $\{w_{11},w_{21},v_3,v_4,v_5\}$ and $\{w_{12},w_{22},v_3,v_4,v_5\}$ induce two $K_5$'s in $R''$. For each $i \in [5]$, note that $\phi^{R''}(W_i)=\phi^R(U_i)=2$ since all vertices receive the same color by Claim \ref{mono}.  Since there is no monochromatic triangle,  the edge colorings of $w_{11},w_{21},v_3,v_4,v_5$ and $w_{12},w_{22},v_3,v_4,v_5$ must be pentagonlike colorings. By symmetry, assume $\phi^{R''}(v_3v_4)=\phi^{R''}(v_3v_5)=1$, which implies that $\phi^{R''}(v_3w_{11})=\phi^{R''}(v_3w_{21})=2$ and $\phi^{R''}(v_3w_{12})=\phi^{R''}(v_3w_{22})=2$.  Since $\phi^{R''}(w_{11}w_{21})$ and $\phi^{R''}(w_{12}w_{22})$ receive different colors, we can easily get a monochromatic triangle, a contradiction. Therefore, for all distinct $i,j\in[5]$, $\phi^{R''}(W_i,W_{j})=1$ or $\phi^{R''}(W_i,W_{j})=2$.

Now, apply $(\ref{4.8})$ again, $R''[W_1,\ldots,W_5]$ must contain a $K_5$ as a subgraph, then the edge coloring of this $K_5$ must be a pentagonlike coloring. The assertion follows.\hfill$\Box$

\medskip
 For each $i\in[6]$,
let $${X_i'} = \sqcup_{\ell \in {U_i}} {{V_\ell}}.$$ Then $V( G ) = \sqcup_{i = 1}^6 {{X_i'}} $. Recall that $\sum\limits_{i \in [6]} {e( {R[ {{U_i}} ]} )}  \le {\gamma ^{\frac{1}{3}}}{m^2}$, $| {{U_i}} | = \frac{m}{6} \pm {\gamma ^{\frac{1}{4}}}m$ from $(\ref{4.5})$ and $|V_i|=\frac{n}{m}$ for $1\le i\le m$, so we obtain that for each $i\in[6]$, $$| {{X_i'}} | = \frac{n}{6} \pm {\gamma ^{\frac{1}{4}}}n,$$  and

\begin{equation}\label{4.9}
\sum\limits_{i \in [ 6 ]} {e( {G[ {{X_i'}} ]} ) \le \sum\limits_{i \in [ 6 ]} {e( {R[ {{U_i}} ]} )} {{\left( {\frac{n}{m}} \right)}^2} + \varepsilon {n^2} + \frac{{{n^2}}}{m} + 2\gamma {n^2} \le 2{\gamma ^{\frac{1}{3}}}{n^2}}.
\end{equation}

%We will consider the following modified partition $\sqcup_{i=1}^6  X_i$ of $V(G)$ to ensure the big crossing degree with respect to the partition $\sqcup_{i=1}^6  X_i$.

\begin{claim}\label{b-c-d}
There exists a partition $\sqcup_{i = 1}^6 {{X_i}} $ of $V(G)$ such that the following hold:

\smallskip
(i)~For each $i\in[6]$, $| {{X_i}} | = \frac{n}{6} \pm 2\gamma ^{\frac{1}{4}}n$, $||X_i|-|X_i'||\leq 36\gamma ^{\frac{1}{3}}n$, and $|X_i\Delta X_i'|\leq 2\gamma ^{\frac{1}{3}}n^2$.

\smallskip
(ii)~${\delta ^{cr}}( {G[  \sqcup_{i=1}^6{{X_i} }  ]} ) \ge \frac{n}{{18}}$.
\end{claim}
{\bf Proof.}~For each $i\in[6]$ and $v\in X_i'$, if $\deg(v,X_j')\leq\frac{n}{18}$ for some $j\neq i$, then move $v$ to $X_j'$. We repeat this operation until no such vertex exists. Let $X_i$ be the resulting set. We first show that $X_i$'s are well-defined. Recall that $\delta(G)\geq \frac{5n}{6}$, and so if there exist distinct $i,j  \in [6]$ and $v\in X_i'$ with $\deg(v,X_j')\leq\frac{n}{18}$, then we see that
\[
\deg( {v,{X_{i}'}} ) \ge \delta ( G ) - \frac{n}{{18}} - 4\mathop {\max }\limits_{\ell \in [ 6 ]} | {{X_\ell'}} | \ge \frac{n}{9}.
\]
Thus, after moving $v$ from $X_i'$ to $X_j'$, the number of inner edges decreasing by at least $\frac{n}{18}$. Hence, from $(\ref{4.9})$, the process will stop after moving at most $2\gamma ^{\frac{1}{3}}n^2/(\frac{n}{18})=36\gamma ^{\frac{1}{3}}n$ vertices. Therefore, we obtain (i). Moreover, (ii) holds by definition.\hfill$\Box$

\medskip
{\bf ($P_1$):} Note that Claim \ref{b-c-d} (i) implies ($P_1$).

\medskip
{\bf ($P_2$):} Recall $\phi^R(U_6)=1$ from Claim \ref{1-v-c} (i), so Claim \ref{mono} implies $\alpha ( {{G_1}[ {{X_6'}} ]} )\le\sqrt \delta  n$. Consequently, by Claim \ref{b-c-d}  (i), we have
\begin{equation}\label{4.10}
  \alpha ( {{G_1}[ {{X_6}} ]} ) \le \alpha ( {{G_1}[ {{X_6'}} ]} ) + | {| {{X_6}} | - | {{X_6'}} |} | \le \sqrt \delta  n + 36{\gamma ^{\frac{1}{3}}}n \le {\gamma ^{\frac{1}{4}}}n.
\end{equation}
Thus ($P_2$) holds.

\medskip
{\bf The first part of ($P_8$):} Similar to $(\ref{4.10})$, we obtain $\alpha ( {{G_2}[ {{X_i}} ]} ) \le {\gamma ^{\frac{1}{4}}}n$ for each $i \in [5]$,  proving the first part of ($P_8$) (the proof of the second part of ($P_8$) will be given in the end of the lemma).

Note that $G_1$ is $K_3$-free, then it follows by $(\ref{4.10})$ that for each $v\in V(G)$, we have

\begin{equation}\label{4.11}
  \deg_{G_1}(v,X_6)\leq \alpha(G_1[X_6])\leq \gamma ^{\frac{1}{4}}n.
\end{equation}

By Claim \ref{b-c-d} (i) and $\sum\limits_{i \in [ 6 ]} {e( {G[ {{X_i'}} ]} ) \le 2{\gamma ^{\frac{1}{3}}}{n^2}}$ from $(\ref{4.9})$, we have

\[
\sum\limits_{i \in [6]} {e( {G[ {{X_i}} ]} )}  \le \sum\limits_{i \in [6]} {( {e( {G[ {{X_i'}} ]} ) + | {| {{X_i}} | - | {{X_i'}} |} |n} )}  \le {\gamma ^{\frac{1}{4}}}{n^2}.
\]
Since $e(G)\geq\frac{5}{12}n^2$ and $e(T_{n,6}) \le \frac{5}{{12}}{n^2}$, the number of non-edges between distinct parts
\[
e( {\overline G [ {{X_1}, \ldots ,{X_6}} ]} ) \le \frac{5}{{12}}{n^2} - \left( {e( G ) - \sum\limits_{i \in [6]} {e( {G[ {{X_i}} ]} )} } \right) \le \sum\limits_{i \in [6]} {e( {G[ {{X_i}} ]} )}  \le {\gamma ^{\frac{1}{4}}}{n^2}.
\]
Define $J$ by the set of vertices with missing crossing degree at least $\gamma ^{\frac{1}{8}}n$, i.e.,
\[
J = \bigcup_{i \in [ 6 ]} \left\{ {v \in  {{X_i}} :\deg(v,V(G)\setminus X_i) \le |V(G)\setminus X_i|  - {\gamma ^{\frac{1}{8}}}n} \right\}.
 \]
From the definition, we know that
\begin{equation}\label{size-J}
  | J | \le \frac{{2e( {\overline G [ {{X_1}, \ldots ,{X_6}} ]} )}}{{{\gamma ^{\frac{1}{8}}}n}}
  = 2{\gamma ^{\frac{1}{8}}}n.
\end{equation}

For convenience, for $A,B\subseteq V$ and $k\in[2]$, denote $\deg_k(v,A)$, $N_k(v,A)$  and $e_{k}(A,B)$ for $\deg_{G_k}(v,A)$, $N_{G_k}(v,A)$ and $e_{G_k}(A,B)$, respectively.
Recall that ${\phi ^{R''}} := { {{\phi ^R}} |_{R''}}.$

\medskip

\begin{proposition}\label{k-4}
If $\phi^{R''}(W_i,W_{j})=k$ for $i,j\in[6]$ and $k\in [2]$, then for any $C_i\subseteq X_i$, $C_j\subseteq X_j$ with $|C_i|, |C_j| \ge 10\gamma ^{\frac{1}{24}}n$, we have $e_{k}(C_i,C_j) \ge \frac{|C_i|^2|C_j|^2}{64n^2}$. %Furthermore, if $\phi^{R''}(W_i,W_{j})=2$ for $i,j\in[5]$, then there exists $\{u_i^1,u_i^2,u_j^1,u_j^2\}$ which induces a $K_{4}$ of color $2$ with $u_i^1,u_i^2 \in X_i$ and $u_j^1,u_j^2 \in X_j$.
\end{proposition}
{\bf Proof.}~ Recall that ${X_i'} = \sqcup_{\ell\in {U_i}} {{V_{\ell}}}$ and $||X_i|-|X_i'||\leq 36\gamma ^{\frac{1}{3}}n$ for all $i\in [6]$ from Claim \ref{b-c-d} (i). Let $C_i'=C_i\cap X_i'$ for $i\in[6]$. Thus we have  $|C_i'|\ge|C_i|-36\gamma ^{\frac{1}{3}}n>|C_i|/2$. For $i\in[6]$, define
$$ Y_i=\left\{\ell: \ell\in U_i \; \text{such that} \; |V_\ell \cap C_i'|\ge \frac12|C_i'| /\left(\frac m6\right) \right\}.$$
Note that $|V_\ell|=\frac{n}{m}$ for all $\ell\in [m]$ and $| {{U_i}} | = \frac{m}{6} \pm {\gamma ^{\frac{1}{4}}}m$ from (\ref{4.5}), it follows that
$$|Y_i|\cdot  \frac nm +\left(|U_i|-|Y_i|\right)\cdot\frac12|C_i'| /\left(\frac m6\right)\ge |C_i'|,$$
implying that
$
|Y_i|(\frac nm-\frac{|C_i'|/2}{m/6})\ge|C_i'|(1-\frac{|U_i|/2}{m/6})\ge|C_i'|/3,
$ and so  we have
$|Y_i|\ge \frac {m|C_i'|}{3n} \ge \frac {m|C_i|}{6n}.$

Recall $R''=R'[ \sqcup_{i=1}^6{{W_i} }]$, where $W_i = {U_i}\backslash ({R_0} \sqcup R_1 )$ for each $i\in [6]$. Then, $$|R''\Delta R'|\le6\cdot|R_0\sqcup R_1|\cdot m\le18\gamma ^{\frac{1}{6}}m^2.$$
 Note that $R'$ is the graph obtained from $R$ by deleting all edges of weight at most $\frac{1}{2}+\gamma$. Recall that $|R\Delta R'|= e(R)-e(R')\le\gamma ^{\frac{1}{4}}m^2$ from $(\ref{4.4})$, and $|R\Delta T_{m,6}|\le\gamma ^{\frac{1}{3}}m^2$ from $(\ref{4.1})$, so we have that $$|R''\Delta T_{m,6}|\le18\gamma ^{\frac{1}{6}}m^2+\gamma ^{\frac{1}{4}}m^2+\gamma ^{\frac{1}{3}}m^2\le\gamma ^{\frac{1}{7}}m^2.$$

For distinct $i,j\in[6]$, $\phi^{R''}(W_i,W_{j})=k$ from the assumption and each edge of $R''$ has weight at least $\frac{1}{2}+\gamma$,
and for all $\ell_i\in Y_i$, $\ell_j\in Y_j$ with $(V_{\ell_i},V_{\ell_j})$ is $\varepsilon$-regular pair, we have
$|V_{\ell_i}\cap C_i'|\ge\frac12|C_i'| /\left(\frac m6\right)>\varepsilon \frac nm$ and
$|V_{\ell_j}\cap C_j'|\ge \frac12|C_j'| /\left(\frac m6\right)>\varepsilon \frac nm$, it follows that
$$d(V_{\ell_i}\cap C_i',V_{\ell_j}\cap C_j') \ge d(V_{\ell_i},V_{\ell_j})-\varepsilon.$$
Note that $||U_i|-|W_i||,||U_j|-|W_j||\le|R_0\sqcup R_1|<3\gamma ^{\frac{1}{6}}m$, so for each $Y_i$, there are at least $(|Y_i|-3\gamma ^{\frac{1}{6}}m)$ vertices belong to $W_i$.
Since $|Y_i|\ge \frac {m|C_i|}{6n}\gg\gamma ^{\frac{1}{6}}m$ and $|C_i'|\ge|C_i|/2$, we have that
\begin{align*}
e_{k}(C_i',C_j')&\ge\bigg(\left(|Y_i|-3\gamma ^{\frac{1}{6}}m\right)\left(|Y_j|-3\gamma ^{\frac{1}{6}}m\right)-\gamma ^{\frac{1}{7   }}m^2\bigg)
 \left(\frac12+\gamma-\varepsilon\right)\frac{|C_i'|/2}{m/6}\cdot\frac{|C_j'|/2}{m/6}
 \\&\ge\frac{|Y_i||Y_j|}{2}\cdot\frac{36|C_i||C_j|}{32m^2}
 \\&\ge \frac{|C_i|^2|C_j|^2}{64n^2}.
\end{align*}
The assertion follows.\hfill$\Box$

\medskip
For $i\in[5]$, define $J_i'=\{v\in X_i:\deg_1(v,X_{i+1})\ge\gamma^{\frac{1}{8}}n\}.$
%recall that $\phi^{R''}(W_i,W_{i+1})=2$ from Claim \ref{1-e-c},
%where the summations of the subscripts are taken modular $5$,
%and we shall show that $|J_i'|$ is small as the following Observation.

\begin{proposition}\label{J'}
For $i\in[5]$, we have $|J_i'|\le\gamma^{\frac{1}{57}}n.$
\end{proposition}
{\bf Proof.}~ For $i\in[6]$, recall that $W_i=U_i\setminus (R_0\sqcup R_1)$, and $|R_0\sqcup R_1|\le 3\gamma^{\frac{1}{6}}m$, and $|W_i|=\frac{m}{6}\pm\gamma^{\frac{1}{7}}m$ from (\ref{4.7}). Since $\phi^{R''}(W_1,W_2)=2$ from Claim \ref{1-v-c} (ii), $\phi^{R''}(\ell_1\ell_2)=\phi^{R}(\ell_1\ell_2)=2$ for each $\ell_1\ell_2\in E(R'')$ with $\ell_1\in W_1$ and $\ell_2\in W_2$. Thus, $d_1(V_{\ell_1},V_{\ell_2})<\gamma$ and $d_2(V_{\ell_1},V_{\ell_2})\ge\gamma$,  implying that
$e_1(\sqcup_{\ell_1\in W_1}V_{\ell_1},\sqcup_{\ell_2\in W_2}V_{\ell_2})\le |W_1||W_2|\cdot\gamma\left(\frac{n}{m}\right)^2\le \left(\frac{m}{6}+\gamma^{\frac{1}{7}}m\right)^2\cdot\gamma\left(\frac{n}{m}\right)^2\le\frac{1}{2}\gamma n^2.$ Therefore,
\begin{align}\label{J_1'}
  e_1(\sqcup_{\ell_1\in U_1}V_{\ell_1},\sqcup_{\ell_2\in U_2}V_{\ell_2}) & \le e_1(\sqcup_{\ell_1\in W_1}V_{\ell_1},\sqcup_{\ell_2\in W_2}V_{\ell_2})+\sum_{i=1}^2\left((|U_i|-|W_i|)\cdot\frac{n}{m}\cdot n\right) \nonumber\\
   & \le \frac{1}{2}\gamma n^2+6\gamma^{\frac{1}{6}}n^2\nonumber
   \\&\le\gamma^{\frac{1}{7}}n^2.
\end{align}

For $i\in[6]$, recall that $X_i'=\sqcup_{\ell\in U_i}V_\ell$, and $||X_i'|-|X_i||\le36\gamma^{\frac{1}{3}}n$ from Claim \ref{b-c-d}. Thus,
\begin{align*}
  |J_1'|\gamma^{\frac{1}{8}}n  \le e_1(X_1,X_2)&\le e_1(X_1',X_2')+ \sum_{i=1}^2\left((|X_i|-|X_i'|)\cdot n\right)\\
  & \le e_1(X_1',X_2')+ 72\gamma^{\frac{1}{3}}n^2\\&\le e_1(\sqcup_{\ell_1\in U_1}V_{\ell_1},\sqcup_{\ell_2\in U_2}V_{\ell_2})+ 72\gamma^{\frac{1}{3}}n^2%\\  &\le\gamma^{\frac{1}{7}}n^2+72\gamma^{\frac{1}{3}}n^2
  \\&\le 2\gamma^{\frac{1}{7}}n^2.
\end{align*}
Therefore, $|J_1'|\le\gamma^{\frac{1}{57}}n.$ Similarly, we have that $|J_i'|\le\gamma^{\frac{1}{57}}n$ for $i\in[5]$ as desired.\hfill$\Box$
\begin{claim}\label{v1-1}
For any $k\in[5]$ and $x\in X_k$,
$\min\{\deg_2(x,X_{k+1}),\deg_2(x,X_{k+2}\})\}\leq 11\gamma ^\frac{1}{58}n$ \;{and}\;
$\min\{\deg_2(x,X_{k+3}),\deg_2(x,X_{k+4})\}\leq 11\gamma ^\frac{1}{58}n$.
\end{claim}
{\bf Proof.}~By symmetry, we consider $k=1$. We only prove the first inequality since the second is similar. For $i\in\{2,3\}$, suppose to the contrary that $\deg_2(x,X_i)\geq 11\gamma ^\frac{1}{58}n$ for some  $x\in X_1$. Recall that $|J|\le2{\gamma ^{\frac{1}{8}}}n$, and  for $i \in [5]$, $|J_i'|\le\gamma^{\frac{1}{57}}n$ from Proposition \ref{J'} and $\alpha ( {{G_2}[ {{X_i}} ]} ) \le {\gamma ^{\frac{1}{4}}}n$ (the first part of ($P_8$)). Thus $\deg_2(x,X_i)\ge |J|+|J_i'|+\alpha(G_2[X_i])$, and so there exists an edge $u_1^2u_2^2$ of color $2$ with $\{u_1^2,u_2^2\}\subseteq N_2(x,X_2)\setminus (J\cup J_2')$.

Let $Y=N_2(x,X_3)$. Since $u_1^2,u_2^2\notin J\cup J_2'$, we have that $\deg(u_1^2,X_3),\deg(u_2^2,X_3)\ge |X_3|-\gamma ^\frac{1}{8}n$ and
$\deg_1(u_1^2,X_3),\deg_1(u_2^2,X_3)\le \gamma ^\frac{1}{8}n$. Thus,
\begin{align*}
  |N_2(\{x,u_1^2,u_2^2\},Y)|& \ge \deg_2(u_1^2,Y)+\deg_2(u_2^2,Y)-|Y| \\
   & =\deg(u_1^2,Y)+\deg(u_2^2,Y)-\deg_1(u_1^2,Y)-\deg_1(u_2^2,Y)-|Y|\\
   &\ge (|Y|-\gamma ^\frac{1}{8}n)+(|Y|-\gamma ^\frac{1}{8}n)-\gamma ^\frac{1}{8}n-\gamma ^\frac{1}{8}n-|Y|\\
   &=|Y|-4\gamma ^\frac{1}{8}n
   \\&\ge\alpha(G_2[X_3])+|J|.
\end{align*}
Therefore, there exists an edge $u_1^3u_2^3$ of color $2$ with $\{u_1^3,u_2^3\}\subseteq N_2(\{x,u_1^2,u_2^2\},Y)\setminus J$. %Since $u_3^1,u_3^2\notin J$, $d(u_3^1,X_6),d(u_3^2,X_6)\ge |X_6|-\gamma ^\frac{1}{8}n$. Together with

Since for each $i\in\{2,3\}$ and $j\in[2]$, $u_j^i\notin J$ and $\deg_{1}(u_j^i,X_6)\leq \gamma ^{\frac{1}{4}}n$ from $(\ref{4.11})$, we have
\[
 {{\deg_2}( {u_j^i,{X_6}} )}={{\deg}( {u_j^i,{X_6}} )}-{{\deg_1}( {u_j^i,{X_6}} )} \ge (| {{X_6}} | - {\gamma ^{\frac{1}{8}}}n) - \gamma ^{\frac{1}{4}}n \ge | {{X_6}} | - 2{\gamma ^{\frac{1}{8}}}n.
\]
Note that ${\delta ^{cr}}( {G[  \sqcup_{i=1}^6{{X_i} }  ]} ) \ge \frac{n}{{18}}$ from Claim \ref{b-c-d} (ii), we obtain that
\[
 {{\deg_2}( {x,{X_6}} )}={{\deg}( {x,{X_6}} )}-{{\deg_1}( {x,{X_6}} )}  \ge {\delta ^{cr}}( {G[ \sqcup_{i=1}^6{{X_i} }  ]} ) - \gamma ^{\frac{1}{4}}n \ge \frac{n}{{18}} - {\gamma ^{\frac{1}{4}}}n.
\]
Consequently,
$
| {{N_2}( {\{ {x,u_1^2,u_2^2,u_1^3,u_2^3} \},{X_6}} )} | \ge 4( {| {{X_6}} | - 2{\gamma ^{\frac{1}{8}}}n} ) + \frac{n}{{18}} - {\gamma ^{\frac{1}{4}}}n - 4| {{X_6}} | \ge \frac{n}{{20}},
$
showing that $K_6\subseteq G_2$, a contradiction.\hfill$\Box$
\begin{claim}\label{v1-2}
For any $k\in[5]$  and $x\in X_k$,
$\deg_2(x,X_{k+2}),\deg_2(x,X_{k-2})\leq 11\gamma ^\frac{1}{58}n$.
\end{claim}
{\bf Proof.} %~ Suppose to the contrary that $d_2(v,X_3)> 5\gamma ^\frac{1}{50}n$ or $d_2(v,X_4)> 5\gamma ^\frac{1}{50}n$ for some vertex $v\in X_1$.
By symmetry, we consider $k=1$. First, suppose  that $\deg_2(x,X_3)> 11\gamma ^\frac{1}{58}n$ and $\deg_2(x,X_4)> 11\gamma ^\frac{1}{58}n$.
Then, by a similar argument as Claim \ref{v1-1}, there exists $\{x,u_1^3,u_2^3,u_1^4,u_2^4,u\}$ which forms a $K_6$ of color $2$ with $u_j^i \in N_2(x,X_i)\setminus J$ for each  $i\in[3,4], j\in[2]$ and $u\in X_6$, a contradiction.

Now, by symmetry, suppose that $\deg_2(x,X_3)> 11\gamma ^\frac{1}{58}n$ but $\deg_2(x,X_4)\leq 11\gamma ^\frac{1}{58}n$.
Then, it follows by Claim \ref{v1-1} that  $\deg_2(x,X_2)\leq 11\gamma ^\frac{1}{58}n$. Thus, Claim \ref{b-c-d} (ii) implies that
\[
 {{\deg_1}( {x,{X_2}} )}={{\deg}( {x,{X_2}} )}-{{\deg_2}( {x,{X_2}} )}  \ge {\delta ^{cr}}( {G[\sqcup_{i=1}^6{{X_i} } ]} ) - 11\gamma ^\frac{1}{58}n > \frac{n}{{19}}\gg 10\gamma ^\frac{1}{24}n,
\]
and similarly,
$
 {{\deg_1}( {x,{X_4}} )} > \frac{n}{{19}}\gg 10\gamma ^\frac{1}{24}n.
$

Note that $\phi^{R''}(W_2,W_4)=1$ from Claim \ref{1-v-c} (ii) and $\deg_1(x,X_2), \deg_1(x,X_4)\ge 10\gamma ^\frac{1}{24}n$,  so there exists an edge $u_2u_4\in E_{1}(N_1(x,X_2),N_1(x,X_4))$ such that $\{x,u_2,u_4\}$ induces a $K_3$ in $G_1$ from Proposition \ref{k-4} by noting $e_{1}(N_1(x,X_2),N_1(x,X_4))\ge\frac{(10\gamma ^\frac{1}{24}n)^4}{64n^2}\ge\gamma ^\frac{1}{6}n^2>0$, which is again a contradiction.
\hfill$\Box$

\medskip

Next, we will show that ($P_3$) and ($P_4$).

\medskip
{\bf ($P_3$):} By symmetry, it is enough to show $\min \{{\deg_{1}(v,X_1), \deg_{1}(v,X_{3})\} \leq \gamma^\frac{1}{59}}n$ for each $v \in X_6$. Recall that $J$ is the set of vertices with a large missing crossing degree and $|J|\le 2\gamma^{\frac18}n$. Suppose to the contrary that $\deg_1(v,X_1),\deg_1(v,X_3)>\gamma ^{\frac{1}{{59}}}n>|J|$ for some vertex $v\in X_6$. From Claim \ref{v1-2}, $\deg_2(w,X_3)\leq 11\gamma ^{\frac{1}{{58}}}n$ for each $w\in X_1$, it follows that for each vertex $w\in N_1(v,X_1)\setminus J$,
$$\deg_1(w,X_3)=\deg(w,X_3)-\deg_2(w,X_3)\geq |X_3|-\gamma ^{\frac{1}{{8}}}n-11\gamma ^{\frac{1}{{58}}}n.$$
Thus,
$
| {{N_1}( {\{ {w,v} \},{X_3}} )} | \ge (| {{X_3}} | - {\gamma ^{\frac{1}{8}}}n - 11{\gamma ^{\frac{1}{{58}}}}n) + {\gamma ^{\frac{1}{{59}}}}n - | {{X_3}} | > {\gamma ^{\frac{1}{{58}}}}n > 0,
$
which implies that there is a $K_3$ in $G_1$, a contradiction.

\medskip
{\bf ($P_4$):} Suppose to the contrary that there exists a vertex $v\in X_6$ such that
$$\min\limits_{i\in [5]} \{{\deg_1(v,X_i\sqcup X_{i+1})\} > \gamma^\frac{1}{60}}n.$$
Then $\deg_1(v,X_1\sqcup X_2)>\gamma ^\frac{1}{60}n$, and so either $\deg_1(v,X_1)>\frac{1}{2}\gamma ^\frac{1}{60}n$ or $\deg_1(v,X_2)>\frac{1}{2}\gamma ^\frac{1}{60}n$. By symmetry, we assume $\deg_1(v,X_1)>\frac{1}{2}\gamma ^\frac{1}{60}n>\gamma ^\frac{1}{59}n$. It follows by ($P_3$) that $\deg_1(v,X_3)\leq \gamma ^\frac{1}{59}n$ and $\deg_1(v,X_4)\leq \gamma ^\frac{1}{59}n$, which implies that $\deg_1(v,X_3\sqcup X_4)\leq 2\gamma ^\frac{1}{59}n <\gamma ^\frac{1}{60}n$, a contradiction.

\medskip
In order to show ($P_5$), we first show the following claims.
\begin{claim}\label{v1-4}
For any $k\in[5]$ and $x\in X_k$, we have
$
\text{$\deg_1(x,X_{k+1}),\deg_1(x,X_{k-1})\leq 11\gamma ^\frac{1}{24}n$}.
$
\end{claim}
{\bf Proof.}~By symmetry, we consider $k=1$. We only prove the first inequality since the second is similar. On  contrary, suppose $\deg_1(x,X_2)> 11\gamma ^\frac{1}{24}n$ for some $x\in X_1$. It follows by Claim \ref{v1-2} and Claim \ref{b-c-d} (ii) that
\begin{align*}
  {\deg_1}( {x,{X_4}} ) & ={\deg}( {x,{X_4}} )-{\deg_2}( {x,{X_4}} ) \\
   & \ge {\delta ^{cr}}( {G[ \sqcup_{i=1}^6{{X_i} } ]} ) - {\deg_2}( {x,{X_4}} ) \\
   &\ge \frac{n}{{18}} - 11{\gamma ^{\frac{1}{{58}}}}n\\
   &\gg 10\gamma ^\frac{1}{24}n.
\end{align*}

Since $\phi^{R''}(W_2,W_4)=1$ from Claim \ref{1-v-c} (ii), and $\deg_1(x,X_2), \deg_1(x,X_4)\ge 10\gamma ^\frac{1}{24}n$, by Proposition \ref{k-4}, there exists an edge $u_2u_4\in E_{1}(N_1(x,X_2),N_1(x,X_4))$ such that $\{x,u_2,u_4\}$ induces a $K_3$ in $G_1$ since $e_{1}(N_1(x,X_2),N_1(x,X_4))\ge\frac{(10\gamma ^\frac{1}{24}n)^4}{64n^2}\ge\gamma ^\frac{1}{6}n^2>0$, a contradiction.\hfill$\Box$

\medskip
Recall that $J := \cup_{i \in [ 6 ]}\{ {v \in  {{X_i}} :\deg(v,V(G)\setminus X_i) \le |V(G)\setminus X_i|  - {\gamma ^{\frac{1}{8}}}n} \}$  consists of vertices with {\em missing crossing degree} at least $\gamma ^{\frac{1}{8}}n$.

\begin{claim}\label{v1-3}
For any  vertex $x\in X_k$ with $k\in[5]$,
${\deg_1}( {x,{X_k}} ) \le | J |\le2\gamma^{\frac18}n.$
\end{claim}
{\bf Proof.}~By symmetry, we show that for any  vertex $x\in X_1$,
${\deg_1}( {x,{X_1}} ) \le | J |\le2\gamma^{\frac18}n.$
Suppose to the contrary that $\deg_1(x,X_1)>|J|$ for some vertex $x\in X_1$. Thus, there exists $xx'\in E(G_1[X_1])$ with $x' \notin J$. Claim \ref{v1-2} implies that
\[
{\deg_1}( {x',{X_3}} ) = \deg( {x',{X_3}} ) - {\deg_2}( {x',{X_3}} ) \ge ( | {{X_3}} | - {\gamma ^{\frac{1}{8}}}n) - 11{\gamma ^{\frac{1}{{58}}}}n \ge | {{X_3}} | - {12\gamma ^{\frac{1}{{58}}}}n.
\]
It follows by Claim \ref{b-c-d} (ii) and Claim \ref{v1-2} that
\[
{\deg_1}( {x,{X_3}} )={\deg}( {x,{X_3}} )-{\deg_2}( {x,{X_3}} ) \ge {\delta ^{cr}}( {G[ \sqcup_{i=1}^6{{X_i} } ]} ) - 11{\gamma ^{\frac{1}{{58}}}}n \ge \frac{n}{{18}} - 11{\gamma ^{\frac{1}{{58}}}}n.
\]
Thus,
$
| {{N_1}( {\{ {x,x'} \},{X_3}} )} | \ge (| {{X_3}} | - 12{\gamma ^{\frac{1}{{58}}}}n) + \left(\frac{n}{{18}} - 11{\gamma ^{\frac{1}{{58}}}}n\right) - | {{X_3}} | > \frac{n}{{20}}>0,
$
showing that $K_3\subset G_1$, a contradiction. Since $| J | \le 2{\gamma ^{\frac{1}{8}}}n$ from (\ref{size-J}), the claim follows.  \hfill$\Box$

\begin{claim}\label{y-X1}
For each vertex $x\in X_k$ with  $k\in[5]$, $\deg_2(x,X_k)\leq 11 \gamma ^\frac{1}{116}n$.
\end{claim}
{\bf Proof.}~By symmetry, we show that for any  vertex $x\in X_1$, $\deg_2(x,X_1)\leq 11 \gamma ^\frac{1}{116}n$.
Suppose to the contrary that there exists a vertex $x\in X_1$ such that $\deg_2(x,X_1)>11 \gamma ^\frac{1}{116}n$.
Define
\[
{Z_4} = \{ {v \in {X_4}:{\deg_2}( {v,{X_1}} ) \ge {\gamma ^{\frac{1}{{116}}}}n} \},
\]
\[
{Z_5} = \{ {v \in {X_5}:{\deg_1}( {v,{X_1}} ) \ge {\gamma ^{\frac{1}{{116}}}}n} \}, \;\;\text{and}\;\;
{Z_6} = \{ {v \in {X_6}:{\deg_1}( {v,{X_1}} ) \ge {\gamma ^{\frac{1}{{116}}}}n} \}.
\]
Recall that $| {{X_i}} | = \frac{n}{6} \pm 2\gamma ^{\frac{1}{4}}n$ for each $i\in[6]$ from Claim \ref{b-c-d} (i).
 It follows by Claim \ref{v1-2} that
\[
| {{Z_4}} | \le \frac{{e( {{G_2}[ {{X_1},{X_4}} ]} )}}{{{\gamma ^{\frac{1}{{116}}}}n}} \le \frac{{| {{X_1}} | \cdot11{\gamma ^{\frac{1}{{58}}}}n}}{{{\gamma ^{\frac{1}{{116}}}}n}} \le 2{\gamma ^{\frac{1}{{116}}}}n.
\]
Similarly, applying Claim \ref{v1-4}, we have
$
| {{Z_5}} | \le \frac{{e( {{G_1}[ {{X_1},{X_5}} ]} )}}{{{\gamma ^{\frac{1}{{116}}}}n}} \le \frac{{| {{X_1}} | \cdot 11{\gamma ^{\frac{1}{{24}}}}n}}{{{\gamma ^{\frac{1}{{116}}}}n}} \le 2{\gamma ^{\frac{1}{{116}}}}n.
$
Note that for each $v\in V(G)$, $\deg_{1}(v,X_6)\leq \gamma ^{\frac{1}{4}}n$ from $(\ref{4.11})$, therefore,
\[
| {{Z_6}} | \le \frac{{e( {{G_1}[ {{X_1},{X_6}} ]} )}}{{{\gamma ^{\frac{1}{{116}}}}n}} \le \frac{{| {{X_1}} | \cdot {\gamma ^{\frac{1}{{4}}}}n}}{{{\gamma ^{\frac{1}{{116}}}}n}} \le 2{\gamma ^{\frac{1}{{116}}}}n.
\]

Since $\deg_1(x,X_5)\leq 11\gamma ^\frac{1}{24}n$ from Claim \ref{v1-4}, and ${\delta ^{cr}}( {G[ \sqcup_{i=1}^6{{X_i} }]} )\ge \frac{n}{18}$ from Claim \ref{b-c-d} (ii), and $\alpha ( {{G_2}[ {{X_i}} ]} )\le \gamma ^{\frac{1}{4}}n$ for each $i\in[5]$ from ($P_2$), it follows that
\begin{align*}
  {\deg_2}( {x,{X_5}} ) &= {\deg}( {x,{X_5}} )-{\deg_1}( {x,{X_5}} ) \\
   & \ge {\delta ^{cr}}( {G[ \sqcup_{i=1}^6{{X_i} }]} ) - 11\gamma ^\frac{1}{24}n\ge \frac{n}{{20}}
   \gg | J | + | {{Z_5}} | + \alpha( {{G_2}[ {{X_5}} ]} ).
\end{align*}
Thus, we can pick $x_1^5,x_2^5\in N_2(x,X_5)\setminus(J\cup Z_5)$ with $\varphi(x_1^5x_2^5)=2$. Since $x_1^5,x_2^5\notin  J$, we obtain $\deg_{1}(x_1^5,X_6), \deg_{1}(x_2^5,X_6)\leq \gamma ^{\frac{1}{4}}n$ from $(\ref{4.11})$. Thus, for each $i\in[2]$,
\[
{\deg_2}( {{x_i^5},{X_6}} )=\deg( {{x_i^5},{X_6}} )-{\deg_1}( {{x_i^5},{X_6}} ) \ge | {{X_6}} | - {\gamma ^{\frac{1}{8}}}n-{\gamma ^{\frac{1}{4}}}n.
\]
Moreover,
$
{\deg_2}( {x,{X_6}} ) ={\deg}( {x,{X_6}} )-{\deg_1}( {x,{X_6}} )\ge {\delta ^{cr}}( {G[ \sqcup_{i=1}^6{{X_i} } ]} ) - {\gamma ^{\frac{1}{4}}}n \ge \frac{n}{{18}} - {\gamma ^{\frac{1}{4}}}n.
$
Therefore, we obtain that
\[
| {{N_2}( {\{ {x,{x_1^5},{x_2^5}} \},{X_6}} )} | \ge \left( \frac{n}{{18}} - {\gamma ^{\frac{1}{4}}}n\right)+2\left( {| {{X_6}} | - {\gamma ^{\frac{1}{{8}}}}n}-{\gamma ^{\frac{1}{4}}}n \right)  - 2| {{X_6}} | \ge \frac{n}{{20}} \gg | J | + | {{Z_6}} |.
\]
Consequently, we can pick a vertex ${x_6} \in {N_2}( {\{ {x,{x_1^5},{x_2^5}} \},{X_6}} )\backslash ( {J \cup {Z_6}} )$ such that
\begin{align*}
   {\deg_2}( {{x_6},N_2(x,X_1)} )& = {\deg}( {{x_6},N_2(x,X_1)} )-{\deg_1}( {{x_6},N_2(x,X_1)} ) \\
  & \ge \left(| N_2(x,X_1) | - {\gamma ^{\frac{1}{8}}}n\right) - {\gamma ^{\frac{1}{{116}}}}n\\
  & \ge | N_2(x,X_1)| - 2{\gamma ^{\frac{1}{{116}}}}n.
\end{align*}
Similarly,
$
{\deg_2}( {{x_1^5},N_2(x,X_1)} ) ,{\deg_2}( {{x_2^5},N_2(x,X_1)} )
\ge | N_2(x,X_1) | - 2{\gamma ^{\frac{1}{{116}}}}n.
$

Now, by Claim \ref{b-c-d} (ii) and $\deg_2(x,X_4)\leq 11\gamma ^\frac{1}{58}n$ from Claim \ref{v1-2}, we have
\begin{align*}
  {\deg_1}( {x,{X_4}} ) & ={\deg}( {x,{X_4}} )-{\deg_2}( {x,{X_4}} ) \\
   & \ge {\delta ^{cr}}( {G[ \sqcup_{i=1}^6{{X_i} } ]} ) - {\deg_2}( {x,{X_4}} ) \\
   &\ge \frac{n}{{18}} - 11{\gamma ^{\frac{1}{{58}}}}n\\
   &\gg | J | + | {{Z_4}} |.
\end{align*}
Thus, we can pick $x_4\in N_1(x,X_4)\setminus(J\cup Z_4)$ such that
\begin{align*}
   {\deg_1}( {x_4,N_2(x,X_1)} )& = {\deg}( {{x_4},N_2(x,X_1)} )-{\deg_2}( {{x_4},N_2(x,X_1)} ) \\
  & \ge \left(| N_2(x,X_1) | - {\gamma ^{\frac{1}{8}}}n\right) - {\gamma ^{\frac{1}{{116}}}}n\\
  & \ge | N_2(x,X_1) | - 2{\gamma ^{\frac{1}{{116}}}}n.
\end{align*}
Let $Y'= N_2( \{x_1^5,x_2^5,x_6,x\},X_1)\cap {N_1}( {x_4,N_2(x,X_1)} )$, we have
\[
| Y' | \ge 4\left( {| N_2(x,X_1) | - 2{\gamma ^{\frac{1}{{116}}}}n}\right ) - 3| N_2(x,X_1) | \ge 3{\gamma ^{\frac{1}{{116}}}}n \ge \delta n \ge\alpha ( G ).
\]
So there exists $x_1^1x_2^1\in E(G[Y'])$. However, if $\varphi(x_1^1x_2^1)=1$, then $\{x_4,x_1^1,x_2^1\}$ induces a $K_3$ in $G_1$, while if $\varphi(x_1^1x_2^1)=2$, then $\{x,x_6,x_1^1,x_2^1,x_1^5,x_2^5\}$ forms a $K_6$ in $G_2$, a contradiction. \hfill$\Box$

\medskip
From the above claim, which together with ${\deg_1}( {x,{X_1}} ) \le | J |\le2\gamma^{\frac18}n$ for each $x\in X_1$ from Claim \ref{v1-3} implies that $\Delta(G[X_1])\leq 11 \gamma ^\frac{1}{116}n+2\gamma^{\frac18}n\le\gamma ^\frac{1}{117}n$. By symmetry, we obtain that
\begin{align}\label{Dt-Xi}
\text{$\Delta(G[X_i])\leq \gamma ^\frac{1}{117}n$ for each $i\in[5]$.}
\end{align}

To complete ($P_5$), it remains to prove the following claim.
\begin{claim}\label{X6}
$\Delta(G[X_6])\leq \gamma ^\frac{1}{117}n.$
\end{claim}
{\bf Proof.}~Define
$
Z_6'=\{v\in X_6:\deg_1(v,X_2)\geq \gamma ^\frac{1}{116}n\}.
$
Since $\deg_{1}(v,X_6)\leq \gamma ^{\frac{1}{4}}n$ for each $v\in V(G)$ from $(\ref{4.11})$ and $| {{X_2}} | = \frac{n}{6} \pm 2\gamma ^{\frac{1}{4}}n$ from Claim \ref{b-c-d} (i), we have
\[
| {{Z_6'}} | \le \frac{{e( {{G_1}[ {{X_2},{X_6}} ]} )}}{{{\gamma ^{\frac{1}{{116}}}}n}} \le \frac{{| {{X_2}} | \cdot {\gamma ^{\frac{1}{{4}}}}n}}{{{\gamma ^{\frac{1}{{116}}}}n}} \le 2{\gamma ^{\frac{1}{{116}}}}n.
\]

On the contrary, suppose there exists a vertex $v\in X_6$ such that $\deg(v,X_6)>\gamma ^{\frac{1}{{117}}}n$. Then $\deg(v,X_6)>\deg_1(v,X_6)+|J\cup Z_6 \cup Z_6'|$, and so we can take a vertex $u\in N_2(v,X_6)\setminus (J\cup Z_6 \cup Z_6')$. By ($P_4$), without loss of generality, we may assume $\deg_1(v,X_1 \sqcup X_2)\leq \gamma^{\frac{1}{60}}n$. Then for $i\in [2]$,
$$\deg_2(v,X_i)=\deg(v,X_i)-\deg_1(v,X_i)\geq {\delta^{cr}}( {G[\sqcup_{i=1}^6{{X_i} } ]} )-\gamma^{\frac{1}{60}}n\ge \frac{n}{18}-\gamma^{\frac{1}{60}}n\geq \frac{n}{20}.$$
By the definition of $J, Z_6
 $ and $Z_6'$,  we have that for $i\in [2]$, $$\deg_2(u,X_i)=\deg(u,X_i)-\deg   _1(u,X_i)\geq \left(|X_i|-\gamma ^{\frac{1}{{8}}}n\right)-\gamma ^{\frac{1}{{116}}}n.$$
  Therefore, we have that for $i\in [2]$, $$|N_2(\{v,u\},X_i)|\geq \left(|X_i|-\gamma ^{\frac{1}{{8}}}n-\gamma ^{\frac{1}{{116}}}n\right)+\frac{n}{20}-|X_i|> \frac{n}{30}\gg |J|+|J_i'|+\alpha(G_2[X_i]).$$
By a similar argument as Claim \ref{v1-1}, there exists $\{u_1^1,u_2^1,u_1^2,u_2^2\}$ which induces a $K_4$ of color $2$ with $\{u_1^1,u_2^1\}\subseteq N_2(\{v,u\},X_1)\setminus(J\cup J_1')$ and $\{u_1^2,u_2^2\}\subseteq N_2(\{v,u,u_1^1,u_2^1\},X_2)$, which in turn implies that $\{v,u,u_1^1,u_2^1,u_1^2,u_2^2\}$ forms a $K_6$ in $G_2$, a contradiction. \hfill$\Box$

\medskip

{\bf ($P_6$):}
Note that $\delta ( G )\ge 5n/6$ and $|X_i|=\frac{n}{6}\pm2\gamma^{\frac14}n$ for $i\in[6]$ from ($P_1$).  Consequently, we can obtain an almost tight crossing degree as
  \[
{\delta ^{cr}}( {G[  \sqcup_{i=1}^6{{X_i} } ]} ) \ge \delta ( G ) - \mathop {\max }\limits_{i \in [ 6 ]} ( {\Delta \left( {G[ {{X_i}} ]} \right) + 4| {{X_i}} |} ) \ge \frac{n}{6} - {\gamma ^{\frac{1}{{118}}}}n,
\]
proving ($P_6$).

\medskip
{\bf ($P_7$):}
Since for each $v\in V(G)$, $\deg_{1}(v,X_6)\leq \gamma ^{\frac{1}{4}}n$ from $(\ref{4.11}) $, we have that for $i\in[5] $ and $v\in X_i$,
$
{\deg_2}( {v,{X_6}} )= {\deg}( {v,{X_6}} )-{\deg_1}( {v,{X_6}} )\ge {\delta ^{cr}}( {G[ \sqcup_{i=1}^6{{X_i} } ]} ) - {\deg_1}( {v,{X_6}} ) \ge | {{X_6}} | - {\gamma ^{\frac{1}{{119}}}}n,
$
proving ($P_7$).

\medskip
{\bf The second part of ($P_8$):}
From Claim \ref{v1-2} and Claim \ref{v1-4}, we have that for $i\in[5] $ and $j_1\in\{i-2,i+2\}$ and $j_2\in\{i-1,i+1\}$, and for any vertex $v\in X_i$,
$$
\text{$\deg_2(v,X_{j_1})\leq 11\gamma ^\frac{1}{58}n$, \;{ and} \; $\deg_1(v,X_{j_2})\leq 11\gamma ^\frac{1}{58}n$}.
$$
 Thus, by ($P_6$), we have
\[
{\deg_1}( {v,{X_{j_1}}} )={\deg}( {v,{X_{j_1}}} )-{\deg_2}( {v,{X_{j_1}}} ) \ge {\delta ^{cr}}( {G[  \sqcup_{i=1}^6{{X_i} }  ]} ) - {\deg_2}( {v,{X_{j_1}}} ) \ge | {{X_{j_1}}} | - {\gamma ^{\frac{1}{{119}}}}n,
\]
and similarly,
$
{\deg_2}( {v,{X_{j_2}}} ) \ge | {{X_{j_2}}} | - {\gamma ^{\frac{1}{{119}}}}n,
$
proving the second part of ($P_8$), as desired.

\subsection{Proof of Theorem \ref{zhu-2}}\label{3--2}

In this subsection, we will give a proof of Theorem \ref{zhu-2}. We restate  it as following.

\medskip\noindent
{\bf Theorem 1.4} \;
Suppose $\frac{1}{n} \ll \delta \ll 1$. Let $G$ be an $n$-vertex $(K_3,K_6)$-free graph with $\alpha(G)\leq\delta n$. Then
$
e(G)\leq \left(\frac{5}{12}+ \frac{\delta }{2}+2.1025\delta^{2}\right) n^{2}.
$

\medskip

We also need the following Lemma, which will be useful to guarantee a certain minimum degree condition in a dense graph.
\begin{lemma}[Kim, Kim and Liu \cite{kkl}]\label{kkl-c}
Suppose $0<\frac{1}{n}\ll \varepsilon \ll d \leq 1$. Suppose that $G$ is an $n$-vertex graph with $e(G)\geq \frac{1}{2}(d+\varepsilon)n^{2}$. Then $G$ contains an $n'$-vertex subgraph $G'$ with $n'\geq \frac{1}{2}\varepsilon^{1/2}n$ such that $e\left( {G'} \right) \ge \frac{1}{2}({d{{n'}^2} + \varepsilon {n^2} - d\left( {n - n'} \right)})$ and $\delta(G')\ge dn'$.
\end{lemma}
\medskip

 \noindent{\bf Proof of Theorem \ref{zhu-2}.}~Suppose to the contrary that $e(G)> (\frac{5}{12}+ \frac{\delta}{2}+\frac{841}{400}\delta^{2}) n^{2}$. We apply Lemma \ref{kkl-c} to obtain an $n'$-vertex graph $G'$ with $n'\geq \frac{1}{2}\delta^\frac{1}{2}n$, $\delta(G')\geq \frac{5n'}{6}$ and
 $e(G') >\frac{1}{2}[\frac{5}{6}n'^2+(\delta+\frac{841}{200}\delta^2)n^2-\frac{5}{6}(n-n')]$.
Let $\delta'=\frac{\delta n}{n'}$. Note that
$
\delta' \in [\delta,\delta ^\frac{1}{3}]
$
as $\delta=\frac{\delta n}{n}\leq \frac{\delta n}{n'}=\delta'\leq\delta n\cdot\frac{2}{\delta^\frac{1}{2}n}=2\delta^\frac{1}{2}\leq \delta ^\frac{1}{3}.$

Since $1\leq\alpha(G')\leq\alpha(G)=\delta n= \delta' n'$, we have that
\begin{align}\label{4.12}
  e(G') & >\frac{1}{2}\left[\frac{5}{6}n'^2+\left(\delta+\frac{841}{200}\delta^2\right)n^2-\frac{5}{6}(n-n')\right]\nonumber\\
   & =\frac{5}{12}n'^2+ \frac{\delta'}{2}n'^2+\frac{841}{400}\delta'^{2} n'^{2}+\frac{1}{2}(\delta n^2-\delta'n'^2)+\frac{841}{400}(\delta^2n^2-\delta'^2n'^2)-\frac{5}{12}(n-n')\nonumber\\
   &=\frac{5}{12}n'^2+ \frac{\delta'}{2}n'^2+\frac{841}{400}\delta'^{2} n'^{2}+\left(\frac{1}{2}\delta n-\frac{5}{12}\right)(n-n')\nonumber\\
   & \ge \frac{5}{12}n'^2+ \frac{\delta'}{2}n'^2+\frac{841}{400}\delta'^{2} n'^{2}.
\end{align}

Note that $\varphi$ still induces an edge-coloring of $G'$ which is $(K_3,K_6)$-free. Since $\frac{1}{n}\ll\delta\ll\gamma$ and $n'\geq \frac{1}{2}\delta^\frac{1}{2}n$ and $\delta' \in [\delta,\delta ^\frac{1}{3}]$, we apply Lemma \ref{color stability} with $G',\delta',\gamma$ (playing the roles of $G,\delta,\gamma$) to obtain a partition $\sqcup_{i=1}^6  X_i$ of $V(G')$ satisfying the following properties.

\medskip
($P_1$) For $i\in [6]$,  $|X_i|=\frac{n'}{6}\pm 2\gamma^{\frac{1}{4}}n'.$

($P_2$) There exists some part, say $X_6$, such that $\alpha(G_1'[X_6])\leq \gamma^{\frac{1}{4}}n'.$

($P_3$) For each $v\in X_6$ and  for $i\in [5]$, we have $\min \{{\deg_{G_1'}(v,X_i), \deg_{G_1'}(v,X_{i+2})\} \leq \gamma^\frac{1}{59}}n'$.

($P_4$) For each $v\in X_6$, we have $\min\limits_{i\in [5]} \{{\deg_{G_1'}(v,X_i\sqcup X_{i+1})\} \leq \gamma^\frac{1}{60}}n'$.

($P_5$) For  $i\in[6]$, we have $\Delta(G'[X_i])\leq \gamma^{\frac{1}{117}}n'.$

($P_6$) $\delta^{cr}(G'[ \sqcup_{i=1}^6{{X_i} } ])\geq \frac{n'}{6}-\gamma^{\frac{1}{118}}n'.$

($P_7$) For $i\in[5]$ and $v\in X_i$, we have $\deg_{G_2'}(v,X_6)\geq| {{X_6}} | - {\gamma ^{\frac{1}{{119}}}}n'.$

($P_8$) For $i \in [5]$ and for each $v \in X_i$,   $\alpha(G_2'[X_i])\leq \gamma^{\frac{1}{4}}n'$, and $\deg_{G_1'}(v,X_{j_1})\geq |X_{j_1}|-\gamma^{\frac{1}{119}}n'$  where $j_1\in \{i-2,i+2\}$ and $\deg_{G_2'}(v,X_{j_2})\geq |X_{j_2}|-\gamma^{\frac{1}{119}}n'$ where $j_2\in \{i-1,i+1\}$.

\medskip

For $i \in [5]$,  define
\[
{I_i} = \left\{ {v \in {X_6}:{\deg_{G_1'}}( {v,{X_i}} )\ge \frac{n'}{10}, \ \text{and} \;\;  {\deg_{G_1'}}( {v,{X_{i + 1}}} ) \ge \frac{n'}{{10}}} \right\},
\]
 and $I=\sqcup_{i \in [5]} {{I_i}}$. Note that $I_i \cap I_{j}= \emptyset$ for all distinct $i, j\in[5]$ from ($P_3$). %To bound $e(G')$, we need to show the following Claims.

\begin{claim}\label{inde-set}
We have that

\smallskip
(i) For $i\in[5]$, the subgraph $G'[X_i]$ is $K_3$-free.

\smallskip

(ii) For $i \in [5]$, $I_i\sqcup I_{i+1}$ forms an independent set, and so $|I_i|+|I_{i+1}|\leq\alpha(G')\leq \delta'n'$.

\end{claim}
{\bf Proof.}~(i) By symmetry, it suffices to show $G'[X_1]$ is $K_3$-free. Suppose to the contrary that $T=\{u,v,w\}$ induces a $K_3$ in $G'[X_1]$. By ($P_1$) and ($P_8$), we have
\[
|N_{G_1'}(T,X_3)|\geq 3(|X_3|-\gamma ^{\frac{1}{{119}}}n')-2|X_3|\geq\frac{n'}{12}.
\]
Since $G_1'$ is $K_3$-free, $T$ is monochromatic in color $2$. Now, use ($P_1$) and ($P_8$) again, we have
\[
|N_{G_2'}(T,X_2)|\geq 3(|X_2|-\gamma ^{\frac{1}{{119}}}n')-2|X_2|\geq\frac{n'}{12}>\alpha(G_2'[X_2]),
\]
which implies that there exists an edge $u_1v_1$ of color $2$ in $N_{G_2'}(T,X_2)$.

For all $i\in[5]$ and each vertex $x\in X_i$, note that $\deg_{G_2'}(x,X_6)\geq| {{X_6}} | - {\gamma ^{\frac{1}{{119}}}}n'$ from ($P_7$),
thus,
\[
|N_{G_1'}(T\cup\{u_1,v_1\},X_6)|\geq 5(|X_6|-\gamma ^{\frac{1}{{119}}}n')-4|X_6|\geq\frac{n'}{12}>0,
\]
which implies that $K_6\subseteq G_2'$, a contradiction, proving  (i).

\medskip
(ii) By symmetry, it suffices to show that $I_1$ and $I_1\sqcup I_2$ are independent.
On contrary, we first suppose that $uv$ is an edge in $G'[I_1]$, from ($P_1$) and the definition of $I_1$, we have
\[
|N_{G_1'}(\{u,v\},X_1)|\geq \frac{n'}{10}+\frac{n'}{10}-|X_1|\geq\frac{n'}{40}.
\]
Since $G_1'$ is $K_3$-free, $\varphi(uv)=2$. Note that $u,v\in I_1\subseteq X_6$, it follows by ($P_3$) that  for each $i\in\{3,4\}$, $\deg_{G_1'}(u,X_i),\deg_{G_1'}(v,X_i)\leq\gamma ^{\frac{1}{{59}}}n'$.

 Since $\delta^{cr}(G[ \sqcup_{i=1}^6{{X_i} } ])\geq \frac{n'}{6}-\gamma^{\frac{1}{118}}n'$ from ($P_6$), we have that  for each $i\in\{3,4\}$,
\[
|N_{G_2'}(\{u,v\},X_i)|\geq 2({\delta ^{cr}}( {[  \sqcup_{i=1}^6{{X_i} }  ]} )-\gamma ^{\frac{1}{{59}}}n')- |X_i|\geq\frac{n'}{7}\gg\alpha(G_2'[X_i]).
\]
Moreover, note that $\deg_{G_2'}(x,X_{4})\geq |X_{4}|-\gamma^{\frac{1}{119}}n'$ for each vertex $x\in X_3$ from ($P_8$), so we can find $u_j^i\in N_{G_2'}(\{u,v\},X_i)$ for each $i\in\{3,4\}$ and $j\in[2]$ such that $\{u,v,u_1^3,u_2^3,u_1^4,u_2^4\}$ induces a  $K_6$ of color $2$, a contradiction. Thus, $I_i$ is an independent set for $i\in[5]$.

Now we suppose that $uv$ is an edge in $G'[I_1\sqcup I_2]$. Since $I_1$ and $I_2$ are independent sets, we may assume that $u\in I_1$ and $v\in I_2$. From the definition of $I_1$ and $I_2$ and ($P_1$), we have
\[
|N_{G_1'}(\{u,v\},X_2)|\geq \frac{n'}{10}+\frac{n'}{10}-|X_2|\geq\frac{n'}{40}>0.
\]
Note that $G_1'$ is $K_3$-free, then $\varphi(uv)=2$. Since $u,v\in X_6$, by ($P_3$), we have $\deg_{G_1'}(u,X_i)\leq\gamma ^{\frac{1}{{59}}}n'$ and $\deg_{G_1'}(v,X_i)\leq\gamma ^{\frac{1}{{59}}}n'$. Thus, by ($P_6$), we have that for each $i\in\{4,5\}$,
\[
|N_{G_2'}(\{u,v\},X_i)|\geq 2({\delta ^{cr}}( {[  \sqcup_{i=1}^6{{X_i} }  ]} )-\gamma ^{\frac{1}{{59}}}n')- |X_i|\geq\frac{n'}{7}\gg\alpha(G_2'[X_i]).
\]
By a similar argument as above, we can find $u_j^i\in N_{G_2'}(\{u,v\},X_i)$ for each $i\in\{4,5\}$ and $j\in[2]$ such that $\{u,v,u_1^4,u_2^4,u_1^5,u_2^5\}$ induces a  $K_6$ of color $2$, a contradiction, proving that $I_i \sqcup I_{i+1}$ forms an independent set for $i\in [5]$.
\hfill$\Box$

\medskip

For each $i\in[5]$,  Claim \ref{inde-set} (i) implies that $\Delta(G'[X_i])\leq\alpha(G')\leq\delta'n'$ and so $$e(G'[X_i])\leq \frac{1}{2}\delta'n'|X_i|.$$
 Recall that $I=\sqcup_{i \in [5]} {{I_i}}$.

\begin{claim}\label{m-X-6}
For each edge $uv\in E(G'[X_6\setminus I])$, we have $\varphi(uv)=1$.
\end{claim}
{\bf Proof.}~Suppose to the contrary that there exists an edge $uv\in E(G'[X_6\setminus I])$ with $\varphi(uv)=2$. Since $u,v\in X_6\setminus I$, we have that for each $i \in [5]$,
\begin{equation}\label{4.121}
  \text{$\min\{\deg_{G_1'}(u,X_i),\deg_{G_1'}(u,X_{i+1})\}<\frac{n'}{10}$,\;{ and}\; $\min\{\deg_{G_1'}(v,X_i),\deg_{G_1'}(v,X_{i+1})\}<\frac{n'}{10}$}.
\end{equation}
 Note that $u,v\in X_6$, by ($P_3$), we have the following two cases:

\smallskip
(a) One of $ \{u,v\}$, say $u$, such that $\deg_{G_1'}(u,X_{j'})>\gamma^\frac{1}{59}n'$ for at most one $j'\in[5]$. Then $\deg_{G_1'}(u,X_{j})\le\gamma^\frac{1}{59}n'$ for $j\in[5]\setminus \{j'\}$.
Since $v\in X_6$,   from ($P_3$), we have that for $i\in [5]$, $\min \{{\deg_{G_1'}(v,X_i), \deg_{G_1'}(v,X_{i+2})\} \leq \gamma^\frac{1}{59}}n'$. Therefore, there is at most one $i_1\in[5]$ such that $\min\{\deg_{G_1'}(v,X_{i_1}),\deg_{G_1'}(v,X_{{i_1}+1})\}>\gamma^\frac{1}{59}n'$.  Together with (\ref{4.121}), there exists $j_1\in[5]$ such that
$\deg_{G_1'}( u,X_{j_1}),\deg_{G_1'}( u,X_{{j_1}+1})\le \gamma^\frac{1}{59}n'$, and $\deg_{G_1'}( v,X_{j_1})\le \gamma^\frac{1}{59}n'$ and $\deg_{G_1'}( v,X_{{j_1}+1})< \frac{n'}{10}$.

\smallskip
(b) Otherwise, for each $w\in \{u,v\}$,  there exists some $i_1\in [5]$ such that $$\min\{\deg_{G_1'}(w,X_{i_1}),\deg_{G_1'}(w,X_{{i_1}+1})\}>\gamma^\frac{1}{59}n'$$ and $\deg_{G_1'}(w,X_{i})\le\gamma^\frac{1}{59}n'$ for  $i\in[5]\setminus\{i_1,i_1+1\}$.
Together with (\ref{4.121}), there exists some $j_1\in[5]$ such that $\deg_{G_1'}( u,X_{j_1})< \frac{n'}{10}$ and $\deg_{G_1'}( u,X_{{j_1}+1})\le \gamma^\frac{1}{59}n'$, and $\deg_{G_1'}( v,X_{j_1})\le \gamma^\frac{1}{59}n'$ and $\deg_{G_1'}( v,X_{{j_1}+1})< \frac{n'}{10}$.

\medskip
Thus we conclude that there always exists some $j_1\in[5]$ such that $\deg_{G_1'}( u,X_{j_1})< \frac{n'}{10}$ and $\deg_{G_1'}( u,X_{{j_1}+1})\le \gamma^\frac{1}{59}n'$, and $\deg_{G_1'}( v,X_{j_1})\le \gamma^\frac{1}{59}n'$ and $\deg_{G_1'}( v,X_{{j_1}+1})< \frac{n'}{10}$. Note from ($P_6$) that ${\delta ^{cr}}( {G'[ \sqcup_{i=1}^6{{X_i} } ]} )\ge \frac{n'}{6}-\gamma^{\frac1{118}}n'$, so we obtain that
\begin{align*}
 {\deg_{{G_2'}}}( {u,{X_{{j_1}}}} ) &={\deg_{{G'}}}( {u,{X_{{j_1}}}} )-{\deg_{{G_1'}}}( {u,{X_{{j_1}}}} )\ge {\delta ^{cr}}( {G'[ \sqcup_{i=1}^6{{X_i} } ]} ) - \frac{n'}{10} \ge \frac{n'}{20},
\\
{\deg_{{G_2'}}}( {v,{X_{{j_1}}}} ) &={\deg_{{G'}}}( {v,{X_{{j_1}}}} )-{d_{{G_1'}}}( {v,{X_{{j_1}}}} )\ge {\delta ^{cr}}( {G'[ \sqcup_{i=1}^6{{X_i} } ]} ) - \gamma^\frac{1}{59}n' \ge \frac{n'}{7},
\end{align*}
and similarly, we have
$\deg_{G_2'}( {v,{X_{{j_1}+1}}} )\ge\frac{n'}{20} \;\; \text{and} \;\;
\deg_{G_2'}( {u,{X_{{j_1}+1}}} )\ge \frac{n'}{7}.$
Therefore, for each ${j} \in \{j_1,j_1+1\}$,
\[
|N_{G_2'}(\{u,v\},X_{{j}})|\geq \frac{n'}{7}+\frac{n'}{20}-|X_{{j}}|\geq \frac{n'}{40}\gg\alpha(G_2'[X_{j}]).
\]

By a similar argument as Claim \ref{inde-set} (ii), we can find $u_{i}^j\in N_{G_2'}(\{u,v\},X_{j})$ for ${j}\in\{j_1,{j_1+1}\}$ and $i\in[2]$ such that $\{u,v,u_1^{j_1},u_2^{j_1},u_1^{j_1+1},u_2^{j_1+1}\}$ induces a  $K_6$ of color $2$, a contradiction.\hfill$\Box$

\medskip
From Claim \ref{m-X-6}, we know that $\varphi(uv)=1$ for any edge $uv\in E(G'[X_6\setminus I])$. Next, we will consider the colors of other edges in $G'[X_6]$.

%If for any edges $uv\in E(X_6)$ satisfying $\phi(uv)=1$, here $u\in \bigcup\limits_{i \in [5]} {{I_i}} $ and $v\in X_6\setminus(\bigcup\limits_{i \in [5]} {{I_i}} )$, then $e(G')\leq\frac{5}{12}n'^2+\frac{1}{2}\delta'n'^2+2\delta'^2n'^2.$
\medskip

\noindent
{\bf Case 1:}
Each edge $uv\in E( G'[I,X_6\setminus I])$ satisfies $\varphi(uv)=1$.

\medskip

Define $A_i=I_i\sqcup I_{i+1}$ and $B_i=I_{i+2}\sqcup I_{i+3}\sqcup I_{i+4}$ for each $i \in [5]$.

\begin{claim}\label{phi 2}
For  $i \in [5]$, $G'[X_6\setminus B_i]$ is $K_3$-free.
\end{claim}
{\bf Proof.}~By symmetry, it suffices to show $K_3\nsubseteq G'[X_6\setminus B_1]$. Otherwise, suppose that $\{x,y,z\}$ induces a triangle in $G'[X_6\setminus B_1]$. Since $G_1'$ is $K_3$-free, we may assume $\varphi(xy)=2$. Note that each edge $uv\in E(G'[I,X_6\setminus I])$ satisfies $\varphi(uv)=1$ from the assumption and $A_1$ is an independent set from Claim \ref{inde-set} (ii), thus $xy\in E(G'[X_6\setminus I])$, contradicting with Claim \ref{m-X-6}.\hfill$\Box$

\medskip
Since $G'[X_i]$ is $K_3$-free for each $i\in[5]$ from Claim \ref{inde-set} (i), we have

$$e(G'[X_i])\leq \frac{1}{2} |X_i|\delta'n',$$
for all $i \in [5]$ by noting $\alpha(G')\le \delta'n'$.

Since $I_i\sqcup I_{i+1}$ is an independent set for each $i\in[5]$ from Claim \ref{inde-set} (ii), we have that
$$e(G'[B_1])=e(G'[I_{3}\sqcup I_{4}\sqcup I_{5}])\leq |I_3||I_5|,$$
and $e(G'[A_1,B_1])\leq |I_1|(|I_3|+|I_4|)+|I_2|(|I_4|+|I_5|).$

Note that $G'[X_6\setminus B_1]$ is $K_3$-free from Claim \ref{phi 2} and the fact that $\alpha(G')\le \delta'n'$, and so
\begin{align*}%\label{X6-A1}
e(G'[X_6\setminus B_1])\leq\frac{1}{2} (|X_6|-|B_1|)\delta'n'.
\end{align*}
Moreover, since $G_1'$ is  $K_3$-free, and $\varphi(uv)=1$ for any $uv\in E(G'[X_6\setminus I])\cup E(G'[I,X_6\setminus I])$ from the assumption and Claim \ref{m-X-6}, we have that
\begin{align*}%\label{A1--X6-A1B1}
  e(G'[B_1,X_6\setminus (A_1\sqcup B_1)])  =e(G'[I_3\sqcup I_4 \sqcup I_5,X_6\setminus (A_1\sqcup B_1)])\nonumber
  \leq\left(\sum\limits_{i = 3}^5 {|I_i|} \right)\delta'n'.
\end{align*}
Thus, %combining (\ref{X6-A1}) with (\ref{A1--X6-A1B1}), we have that
\begin{align*}
  e(G'[X_6]) & = e(G'[X_6\setminus B_1])+e(G'[B_1])+e(G'[A_1,B_1])+e(G'[B_1,X_6\setminus (A_1\sqcup B_1)])\\
  % & \leq \frac{1}{2}(|X_6|-|I_3|-|I_4|-|I_5|)\delta'n'+|I_3||I_5|+|I_1||I_3|+|I_1||I_4|+|I_2||I_4|+|I_2||I_5|+(|I_3|+|I_4|+|I_5|)\delta'n'\\
   & \leq \frac{1}{2}|X_6|\delta'n'+\left(\sum\limits_{i = 1}^5 {|I_i||I_{i+2}|} \right)+\frac{1}{2}\left(\sum\limits_{i = 3}^5 {|I_i|}\right)\delta'n'\\
  & \mathop  \le \limits^{(*)}  \frac{1}{2}|X_6|\delta'n'+2\delta'^2n'^2,
\end{align*}
where the last inequality holds from the computation by LINGO in the Appendix by noting that $|I_i|+|I_{i+1}|\le \delta'n'$ for $i\in[5]$ from Claim \ref{inde-set} (ii).
Therefore,
\[
e(G')\leq e(G'[X_1,\ldots,X_6])+\sum\limits_{i \in [6]} {e( {G'[{X_i}]} )}\leq \frac{5}{12}n'^2+\frac{1}{2}\delta'n'^2+2\delta'^2n'^2,
\]
which leads to a contradiction from (\ref{4.12}).

\medskip\noindent
{\bf Case 2:} There exists an edge $uv\in E(G'[I,X_6\setminus I])$ satisfying $\varphi(uv)=2$.

\medskip
We first have the following four claims.

\begin{claim}\label{4-1-14}
If $u\in I_\ell$ for some $\ell\in [5]$, then $ {{\deg_{G_1'}}( {v,{X_{\ell + 3}}} )}  \ge \frac{n'}{{10}}$.% Furthermore, we have that $ {{d_{G_2'}}( {u',{X_i}} )}  \ge \frac{n'}{{16}}$ for all $i\in [5]\setminus\{j+3\}$, where the summation of the subscript is taken modular 5.
\end{claim}
{\bf Proof.}~Suppose to the contrary that $ {{\deg_{G_1'}}( {v,{X_{\ell + 3}}} )}  < \frac{n'}{{10}}$. Thus, by ($P_6$), we obtain
\begin{align}\label{d-u-j+3}
{{\deg_{G_2'}}( {v,{X_{\ell + 3}}} )}={{\deg_{G'}}( {v,{X_{\ell + 3}}} )}-{{\deg_{G_1'}}( {v,{X_{\ell + 3}}} )}  \ge {\delta ^{cr}}( {G'[  \sqcup_{i=1}^6{{X_i} }  ]} ) - \frac{n'}{{10}} > \frac{n'}{{16}}.
\end{align}
Since $v\in X_6$, we have $\min \{{\deg_{G_1'}(v,X_{\ell+2}), \deg_{G_1'}(v,X_{\ell+4})\} \leq \gamma^\frac{1}{59}}n'$ from ($P_3$).
By symmetry, we may assume ${\deg_{G_1'}}( {v,{X_{\ell + 2}}} ) \le {\gamma ^{\frac{1}{{59}}}}n'$, which implies that
\begin{align}\label{d-u-j+2}
 {{\deg_{G_2'}}( {v,{X_{\ell + 2}}} )}={{\deg_{G'}}( {v,{X_{\ell + 2}}} )}-{{\deg_{G_1'}}( {v,{X_{\ell + 2}}} )}  \ge {\delta ^{cr}}( {G'[ \sqcup_{i=1}^6{{X_i} } ]} ) - {\gamma ^{\frac{1}{{59}}}}n' \ge \frac{n'}{7}.
\end{align}

Note that $u\in I_\ell\subseteq X_6$, by the definition, we have that
\[
{{\deg_{G_1'}}( {u,{X_\ell}} ) \ge \frac{n'}{{10}}, \ \text{and} \;\; {\deg_{G_1'}}( {u,{X_{\ell + 1}}} ) \ge \frac{n'}{{10}}},
\]
and it follows from ($P_3$) that ${\deg_{G_1'}}( {u,{X_{\ell+2}}} ), {\deg_{G_1'}}( {u,{X_{\ell+3}}} )\le {\gamma ^{\frac{1}{{59}}}}n'$, which implies that $$ {{\deg_{G_2'}}( {u,{X_{\ell+2}}} )} ={{\deg_{G'}}( {u,{X_{\ell+2}}} )}-{{\deg_{G_1'}}( {u,{X_{\ell+2}}} )} \ge {\delta ^{cr}}( {G'[ \sqcup_{i=1}^6{{X_i} }  ]} ) - {\gamma ^{\frac{1}{{59}}}}n' \ge \frac{n'}{7},$$
and similarly, $ {{\deg_{G_2'}}( {u,{X_{\ell+3}}} )}\ge \frac{n'}{7}.$
Therefore,  together with (\ref{d-u-j+2}), we have $$|N_{G_2'}(\{u,v\},X_{\ell+2})|\ge \frac{n'}{7}+\frac{n'}{7}-|X_{\ell+2}|>\frac{n'}{10}\gg\alpha(G_2'[X_{\ell+2}]),$$
and recall (\ref{d-u-j+3}), we have
$$|N_{G_2'}(\{u,v\},X_{\ell+3})|\ge \frac{n'}{16}+\frac{n'}{7}-|X_{\ell+3}|>\frac{n'}{50}\gg\alpha(G_2'[X_{\ell+3}]),$$

By a similar argument as Claim \ref{inde-set} (ii), there exists some $u_i^{j}\in N_{G_2'}(\{u,v\},X_i)$ for each $i\in\{\ell+2,\ell+3\}$ and ${j}\in[2]$ such that $\{u,v,u_{\ell+2}^1,u_{\ell+2}^2,u_{\ell+3}^1,u_{\ell+3}^2\}$ induces a  $K_6$ of color $2$, a contradiction.\hfill$\Box$%, proving the first part of Claim.

%Note that $u'\notin I_{j+2}, I_{j+3}$ and $ {{d_{G_1'}}( {u',{X_{j + 3}}} )} \ge \frac{n'}{{10}}$, then
%$$
%{{d_{G_1'}}( {u',{X_{j+2}}} )}, {{d_{G_1'}}( {u',{X_{j+4}}} )}<\frac{n'}{{10}}.
%$$
%Thus,
%$$ {{d_{G_2'}}( {u',{X_{i}}} )}={{d_{G'}}( {u',{X_{i}}} )}-{{d_{G_1'}}( {u',{X_{i}}} )}  \ge {\delta ^{cr}}( {G'[  \sqcup_{i=1}^6{{X_i} }  ]} ) - \frac{n'}{{10}} \ge \frac{n'}{{16}},$$
%for each ${i}\in \{j+2,j+4\}$.

%It follows from $ {{d_{G_1'}}( {u',{X_{j + 3}}} )}  \ge \frac{n'}{{10}}$ have shown from the first part of this Claim and \textbf{(A4)} that
%$$
%{{d_{G_1'}}( {u',{X_{j}}} )}, {{d_{G_1'}}( {u',{X_{j+1}}} )}\le{\gamma ^{\frac{1}{{25}}}}n',
%$$
%which implies that
%$$ {{d_{G_2'}}( {u',{X_{i}}} )} ={{d_{G'}}( {u',{X_{i}}} )}-{{d_{G_1'}}( {u',{X_{i}}} )} \ge {\delta ^{cr}}( {G'[  \sqcup_{i=1}^6{{X_i} } ]} ) - {\gamma ^{\frac{1}{{25}}}}n'\ge \frac{n'}{{16}},$$
 %for all ${i}\in \{j,j+1\} $, proving the second part of Claim.

\begin{claim}\label{4-1-15}
If $u\in I_\ell$ for some $\ell\in [5]$, then $\{v\} \sqcup I_{\ell+2} \sqcup I_{\ell+3 }$ is an independent set.
\end{claim}
{\bf Proof.}~We first show $N_{G'}(v,I_{\ell+2})=\emptyset$. Otherwise, suppose $w\in N_{G'}(v,I_{\ell+2})$. Note that $\deg_{G_1'}(v,X_{\ell+3})\ge\frac{n'}{10}$ from Claim \ref{4-1-14} as $u\in I_\ell$, and from the definition $\deg_{G_1'}(w,X_{\ell+3})\ge\frac{n'}{10}$ as $w\in I_{\ell+2}$, and so we have that \[|N_{G_1'}(\{v,w\},X_{\ell+3})|\ge \frac{n'}{10}+\frac{n'}{10}-|X_{\ell+3}|>0.\] Since $G_1'$ is $K_3$-free, $\varphi(vw)=2$. Thus, by Claim \ref{4-1-14}, $\deg_{G_1'}(v,X_{(\ell+2)+3})=\deg_{G_1'}(v,X_{\ell})\ge \frac{n'}{10}$. %Since $u\in I_\ell$, we have that $d_{G_1'}(v,X_{\ell+3})\ge \frac{n'}{10}$ from Claim \ref{4-1-14}.
Moreover,  $\deg_{G_1'}(v,X_{\ell+3})\ge\frac{n'}{10}$, and $\deg_{G_1'}(x,X_{\ell})\geq |X_{\ell}|-\gamma^{\frac{1}{119}}n'$ for each vertex $x\in X_{\ell+3}$ from ($P_8$), so we can find  $x_{i}\in N_{G_1'}(v,X_{{i}})$ for  ${i}\in \{\ell,\ell+3\}$ such that $\{v,x_\ell,x_{\ell+3}\}$ induces a $K_3$ of color $1$, a contradiction.

By a similar argument as above, we can show $N_{G'}(v,I_{\ell+3})=\emptyset$. Recall that $I_{\ell+2}\sqcup I_{\ell+3}$ is an independent set from  Claim \ref{inde-set} (ii), so the assertion follows. \hfill$\Box$

\medskip

For each $\ell\in[5]$, define $Z_\ell=\{v\in X_6\setminus I: \varphi(uv)=2 ~for~some~u\in I_\ell\}$.

\begin{claim}\label{4-1-16}
For $\ell \in [5]$, $Z_\ell$ is an independent set. %and $Z_i \cap Z_{i+1}=\emptyset$. In particular,
Furthermore, we have that $Z_\ell \sqcup I_{\ell+2} \sqcup I_{\ell+3}$ is an independent set, and so $|Z_\ell|+|I_{\ell+2}|+|I_{\ell+3}| \le \delta'n'$.
\end{claim}
{\bf Proof.}~By Claim \ref{4-1-15},  it suffices to show $Z_\ell$ is an independent set for $\ell \in [5]$. By symmetry, we show $Z_1$ is an independent set. Otherwise, suppose $xy\in E(G'[Z_1])$. %, and $v\in Z_1\cap Z_2$, there exist $v_i\in Z_i$ such that $\phi(vv_i)=2$ for $i\in[2]$.
Note that $Z_1\subseteq X_6\setminus I$, we have $\varphi(xy)=1$ by Claim \ref{m-X-6}. Since $x,y\in Z_1$, we have that $\deg_{G_1'}(x,X_4), \deg_{G_1'}(y,X_4)\ge \frac{n'}{10}$ from Claim \ref{4-1-14}, implying that $\deg_{G_1'}(\{u,v\},X_4) \ge \frac{2}{10}n'-|X_4|>0$. Thus we get a $K_3$ in $G_1'$, a  contradiction. So $|Z_\ell|+|I_{\ell+2}|+|I_{\ell+3}| \le \delta'n'$ by noting $\alpha(G')\le\delta n=\delta'n'$.
\hfill$\Box$

\begin{claim}\label{4-1-17}
For all $i\neq j\in [5]$, $Z_i\cap Z_j=\emptyset$.
\end{claim}
{\bf Proof.}~By symmetry, it suffices to show $Z_1\cap Z_2=\emptyset$ and $Z_1\cap Z_3=\emptyset$. We first show $Z_1\cap Z_2=\emptyset$. Otherwise, suppose $v\in Z_1\cap Z_2$. Recall the definition of $Z_1,Z_2$ and Claim \ref{4-1-14}, we have $ {{\deg_{G_1'}}( {v,{X_{4}}} )}, {{\deg_{G_1'}}( {v,{X_{5}}} )}  \ge \frac{n'}{{10}}$, implying that $v\in I_4$, which is impossible since $v\in Z_1\cap Z_2\subseteq X_6\setminus I$. It remains to show $Z_1\cap Z_3=\emptyset$. Otherwise, suppose $v\in Z_1\cap Z_3$. Recall the definition of $Z_1,Z_3$ and Claim \ref{4-1-14}, we have $ {{\deg_{G_1'}}( {v,{X_{4}}} )}, {{\deg_{G_1'}}( {v,{X_{3+3}}} )}={{\deg_{G_1'}}( {v,{X_{1}}} )}  \ge \frac{n'}{{10}}$. Note that $\deg_{G_1'}(x,X_1)\geq |X_{4}|-\gamma^{\frac{1}{119}}n'$ for each vertex $x\in X_{4}$ from ($P_8$), so we can find  $x_{i}\in N_{G_1'}(v,X_{{i}})$ for  ${i}\in \{1,4\}$ such that $\{v,x_1,x_{4}\}$ induces a $K_3$ of color $1$, a contradiction.\hfill$\Box$
\medskip

Now, we define $Z=\sqcup_{i=1}^5Z_i$.
Since $G'[X_i]$ is $K_3$-free for each $i\in[5]$ from Claim \ref{inde-set} (i), we have that for each $i \in [5]$,
$$e(G'[X_i])\leq \frac{1}{2} |X_i|\delta'n'$$ by noting  $\alpha(G')\le \delta'n'$.

\begin{claim}\label{degree-v}
For each $v\in X_6\setminus I$ and $i\in[5]$,
\begin{align*}
  {\deg_{G_1'}}(v,X_6\setminus I)-{\deg_{G_1'}}(v,Z_i\sqcup Z_{i+1})+{\deg_{G_1'}}(v,I_i\sqcup I_{i+1})\le \delta'n'.
\end{align*}
\end{claim}
{\bf Proof.}~Recall that $\varphi(uw)=1$ for any $uw\in E(G'[X_6\setminus I])$ from Claim \ref{m-X-6}, then $\varphi(e)=1$ for each edge $e$ in $N_{G'_1}(v,X_6\setminus I)$. Furthermore, from the definition of $Z_i$ and $Z_{i+1}$, $\varphi(e)=1$ for each edge $e$ in $N_{G'_1}(v,(X_6\setminus I)\setminus (Z_i\sqcup Z_{i+1})\sqcup (I_i \sqcup I_{i+1}))$. Note that $G_1'$ is $K_3$-free, and $\alpha(G')\le \delta'n'$, the assertion follows since otherwise there exists $xy$ in $N_{G'_1}(v,(X_6\setminus I)\setminus (Z_i\sqcup Z_{i+1})\sqcup (I_i \sqcup I_{i+1}))$ with $\varphi(xy)=1$, and so $\{v,x,y\}$ forms a $K_3$ of color $1$, a contradiction.\hfill$\Box$

\medskip

It follows from Claim \ref{degree-v} that  \begin{align}\label{center fomu}
   \sum\limits_{i = 1}^5\left(\sum\limits_{v \in {X_6}\backslash I} \left[\deg_{G_1'}(v,X_6\setminus I)-\deg_{G_1'}(v,Z_i\sqcup Z_{i+1})+\deg_{G_1'}(v,I_i\sqcup I_{i+1})\right]\right) \le  5(|X_6|-|I|)\delta'n'.
 \end{align}

%From Claim \ref{degree-v}, and $Z_i\cap Z_{i+1}=\emptyset$ for all $i\in[5]$ from Claim \ref{4-1-17} and $I_i\cap I_{i+1}=\emptyset$ for all $i\in[5]$ from ($P_3$), we obtain that
% we have that ${\deg_{G_1'}}(v,Z_i\sqcup Z_{i+1})={\deg_{G_1'}}(v,Z_i)+{\deg_{G_1'}}(v,Z_{i+1})$ and ${\deg_{G_1'}}(v,I_i\sqcup I_{i+1})={\deg_{G_1'}}(v,I_i)+{\deg_{G_1'}}(v,I_{i+1})$, which implies

 %Thus, $e(G_1'[Z_i,X_6\setminus I])\le |Z_i|\delta'n'$ for each $i\in[5]$.

For the first term of (\ref{center fomu}), %since $\sum\limits_{v \in {X_6}\backslash I}{\deg_{G_1'}}(v,X_6\setminus I)=2e(G_1'[X_6\setminus I])=2e(G'[X_6\setminus I])$,
we have $$\sum\limits_{i = 1}^5\left(\sum\limits_{v \in {X_6}\backslash I} \left[\deg_{G_1'}(v,X_6\setminus I)\right]\right)=5\cdot2e(G_1'[X_6\setminus I])=10e(G_1'[X_6\setminus I]).$$

For the second term of (\ref{center fomu}), recall that for each edge $uv\in E(G'[X_6\setminus I])$, $\varphi(uv)=1$ from Claim \ref{m-X-6}, and $G_1'$ is $K_3$-free, and $\alpha(G')\le \delta'n'$, and $Z_i$ is an independent set for each $i\in[5]$ from Claim \ref{4-1-16}. Thus, $e(G_1'[Z_i,X_6\setminus I])=\sum\limits_{v \in Z_i}{\deg_{G_1'}}(v,X_6\setminus I)\le \sum\limits_{v \in Z_i}\delta'n'= |Z_i|\cdot \delta'n'$, which implies that $\sum\limits_{v \in {X_6}\backslash I}{\deg_{G_1'}}(v,Z_i)=e(G_1'[Z_i,X_6\setminus I])\le |Z_i|\cdot \delta'n'$. Hence, by noting $Z_i\cap Z_{i+1}=\emptyset$ for all $i\in[5]$ from Claim \ref{4-1-17}, we obtain that
\begin{align*}
  \sum\limits_{i = 1}^5\left(\sum\limits_{v \in {X_6}\backslash I} \left[{\deg_{G_1'}}(v,Z_i\sqcup Z_{i+1})\right]\right) \le \sum\limits_{i = 1}^5 (|Z_i|+|Z_{i+1}|)\cdot \delta'n'=2|Z|\cdot \delta'n'.
\end{align*}

For the third term of (\ref{center fomu}), since $I_i\cap I_j=\emptyset$ for all $i,j\in[5]$ from ($P_3$), we have
\begin{align*}
  \sum\limits_{i = 1}^5\left(\sum\limits_{v \in {X_6}\backslash I} \left[{\deg_{G_1'}}(v,I_i\sqcup I_{i+1})\right]\right) & =\sum\limits_{i = 1}^5\left(\sum\limits_{v \in {X_6}\backslash I} \left({\deg_{G_1'}}(v,I_i)+{\deg_{G_1'}}(v,I_{i+1})\right)\right)  \\
   & = \sum\limits_{i = 1}^5 \left(e(G_1'[I_i,X_6\setminus I])+e(G_1'[I_{i+1},X_6\setminus I])\right)\\
   &=2e(G_1'[I,X_6\setminus I]).
\end{align*}

Hence, (\ref{center fomu}) implies that
\begin{align}\label{center-f}
  10e(G'[X_6\setminus I])+2e(G_1'[I,X_6\setminus I]) \le 2|Z|\delta'n'+5(|X_6|-|I|)\delta'n'.
\end{align}
By applying Claim \ref{m-X-6}, $e(G_1'[I_i,X_6\setminus I])=\sum\limits_{v \in I_i}{\deg_{G_1'}}(v,X_6\setminus I)\le \sum\limits_{v \in I_i}\delta'n'= |I_i|\cdot \delta'n',$ and so $e(G_1'[I,X_6\setminus I])=\sum\limits_{i = 1}^5e(G_1'[I_i,X_6\setminus I])\le \sum\limits_{i = 1}^5|I_i|\cdot \delta'n'=|I|\cdot \delta'n'.$
Therefore, (\ref{center-f}) implies that
\begin{align*}
  10e(G'[X_6\setminus I])+10e(G_1'[I,X_6\setminus I]) & \le 2|Z|\delta'n'+5(|X_6|-|I|)\delta'n'+8e(G_1'[I,X_6\setminus I])\\
  &\le 2|Z|\delta'n'+5(|X_6|-|I|)\delta'n'+8|I|\delta'n'\\
   & = 5|X_6|\delta'n'+2|Z|\delta'n'+3|I|\delta'n',
\end{align*}
implying
\begin{align}\label{I,X_6I,1}
  e(G'[X_6\setminus I])+e(G_1'[I,X_6\setminus I]) & \le \frac{1}{2}|X_6|\delta'n'+ \frac{1}{5}|Z|\delta'n'+ \frac{3}{10}|I|\delta'n'.
\end{align}

 %by noting that $N_{G_1'}(u,X_6\setminus I)\le\delta'n'$ for each $u\in I$ since $\alpha(G'[X_6\setminus I])\le\alpha(G')\le \delta'n'$.

Recall $Z_i=\{v\in X_6\setminus I: \varphi(uv)=2 ~for~some~u\in I_i\}$ for $i\in[5]$, so we have
\begin{align}\label{I,X_6I,2}
  e(G_2'[I,X_6\setminus I]) & \le \sum\limits_{i = 1}^5 {{|Z_i|}| {{I_i}} |}.
\end{align}
Since $I_i\sqcup I_{i+1}$ is an independent set for each $i\in[5]$ from Claim \ref{inde-set} (ii), we have  that
\begin{align}\label{I}
  e(G'[I])&\le\sum_{i=1}^{5}|I_i||I_{i+2}|.
\end{align}
Therefore, combining with (\ref{I,X_6I,1}), (\ref{I,X_6I,2}) and (\ref{I}), we obtain that
\begin{align*}
  e(G'[X_6]) & =e(G'[X_6\setminus I])+e(G'[I,X_6\setminus I])+e(G'[I])\\
             &=\left(e(G'[X_6\setminus I])+e(G_1'[I,X_6\setminus I])\right)+e(G_2'[I,X_6\setminus I])+e(G'[I])\\
         & \le \frac{1}{2}|X_6|\delta'n'+ \frac{1}{5}|Z|\delta'n'+ \frac{3}{10}|I|\delta'n'+\sum\limits_{i = 1}^5 {{|Z_i|}| {{I_i}} |}+\sum_{i=1}^{5}|I_i||I_{i+2}|\\
         &\mathop  \le \limits^{(**)}  \frac{1}{2}|X_6|\delta'n'+ \frac{841}{400}\delta'^2n'^2,
\end{align*}
where the last inequality holds from the computation by LINGO in the Appendix by noting that $|Z_i|+|I_{i+2}|+|I_{i+3}|\le \delta'n'$ for $i\in[5]$ from Claim \ref{4-1-16}. Consequently,
\[
e(G')\leq e(G'[X_1,\ldots,X_6])+\sum\limits_{i \in [6]} {e( {G'[{X_i}]} )}\leq \frac{5}{12}n'^2+\frac{1}{2}\delta'n'^2+\frac{841}{400}\delta'^2n'^2,
\]
which leads to a contradiction from (\ref{4.12}) again. This completes the proof of Case 2 and hence Theorem \ref{zhu-2}.\hfill$\Box$

\section{Concluding remarks and problems}\label{crp}

In this paper, we make a substantial step toward Conjecture \ref{kkl-2} due to Kim, Kim and Liu \cite{kkl} by showing that for any sufficiently small $\delta>0$, $\rho(3,6,\delta)\le\frac{5}{{12}} + \frac{\delta }{2}+ 2.1025\delta ^2.$
However, there remains a small gap compared with the lower bound $\rho(3,6,\delta)\ge\frac{5}{{12}} + \frac{\delta }{2}+ 2\delta ^2.$ It's worth mentioning that if all edges in $G'[I,X_6\setminus I]$ are in color 1, then  Conjecture \ref{kkl-2} would follow from a slight modification of the proof of Case 2 in Section \ref{3--2}.
% Conjecture \ref{kkl-2} is still open.

We now turn our attention to the next value, $\rho(3,7,\delta).$
Let $d,n$ be integers. Recall that $F(n,d)$ is an $n$-vertex $d$-regular $K_3$-free graph with independence number $d$, which is well defined from \cite{Brandt}.
Suppose that $8$ divides $n$. Let $G$ be a graph obtained from $T_{n,8}$ by placing a copy of $F(\frac{n}{8},d)$, for some $d\in[\delta n-o(n),\delta n]$, in each partite of $T_{n,8}$. It is easy to see that $\alpha(G)\le \delta n$ and $e(G)=\frac{7}{16}n^2+\frac{\delta}{2}n^2+o(n^2)$. Define an edge-coloring $\phi$ of $G$ as follows:

\medskip
(1) $\phi(e)=2$ for all $e\in\cup_{i\in[8]}G[X_i]$;

\smallskip
(2) $\phi(X_i,X_j)=2$ iff $|i-j|\in[2]$ for all $i,j\in[8]$;

\smallskip
(3) all other edges receive color $1$.

\medskip
Then $\phi$ is a $(K_3,K_7)$-free coloring, which implies that $\rho(3,7,\delta)\ge\frac{7}{16}+\frac{\delta}{2}.$

\medskip
We believe the lower bound is tight.

\begin{conjecture}
For sufficiently small $\delta>0$,
 $\rho(3,7,\delta)= \frac{7}{16}+\frac{\delta}{2}.$
\end{conjecture}

\ignore{
We know that the values of $\rho(3,q)$ for $q\ge8$ are wide open.
Let us conclude this section with the following conjecture.
\begin{conjecture}[\cite{B3,kkl}]
For each $q \ge 5$,
\[
\left\{ \begin{array}{ll}
\rho ( {3,2q - 1} ) = \frac{1}{2}\left(1 - \frac{1}{r(3,q) - 1}\right), &\vspace{2mm}\\
\rho ( {3,2q}) = \frac{1}{2}\left(1 - \frac{1}{r(3,q)}\right). &
\end{array} \right.
\]
\end{conjecture}
}

\medskip\noindent
{\bf Acknowledgement:} The authors would like to thank Hong Liu for invaluable discussions and comments, and Ping Hu for carefully reading the manuscript and pointing out a number of details that required correction.

%\medskip

\newpage
\section*{Appendix}
\noindent
All summations of the subscripts are taken modular 5. Define two functions $f$ and $g$ as follows:
\begin{align*}
  f(x_1,\ldots,x_5,y_1,\ldots,y_5) & =\frac{3}{10}\left(\sum\limits_{i = 1}^5 {x_i} \right)+\frac{1}{5}\left(\sum\limits_{i = 1}^5 {y_i} \right) +\left(\sum\limits_{i = 1}^5 {x_iy_{i+2}} \right)+\left(\sum\limits_{i = 1}^5 {x_ix_{i+2}} \right),
\end{align*}
and
\[
g(x_1,\ldots,x_5)=\frac{1}{2}\left(\sum\limits_{i = 3}^5 {x_i} \right)+\left(\sum\limits_{i = 1}^5 {x_ix_{i+2}} \right).
\]
\medskip

The domains of $f$ and $g$ are
\[
D_f=\{(x_1,\ldots,x_5,y_1,\ldots,y_5): x_i, y_i\ge 0 ~\text{and}~ x_i+x_{i+1}+y_{i}\le 1 ~\forall~ i\in [5]\},
\]
and
\[
D_g=\{(x_1,\ldots,x_5): x_i\ge 0~\text{and}~ x_i+x_{i+1}\le 1 ~\forall~ i\in [5]\},
\]
respectively.

\medskip

Applying the LINGO, we have
\begin{equation}\label{5-1}
\mathop {\max }\limits_{{D_f}} f=f(0.45,0.55,0.45,0,0,0,0,0.55,1,0)=2+\frac{41}{400}.
\end{equation}

\begin{equation}\label{5-2}
\mathop {\max }\limits_{{D_g}} g=g(0.5,0.5,0.5,0.5,0.5)=2.
\end{equation}

We first show  $(*)$ in the end of Case 1 of the proof of Theorem \ref{zhu-2}. %and $(**)$.
\medskip

Let $x_i=|I_i|/(\delta'n')$.  From Claim \ref{inde-set} (ii),  for $i\in [5]$, $$|I_i|+|I_{i+1}|\le \delta'n'.$$ Then  for $i\in [5]$, $$x_i+x_{i+1}\le 1,~ x_i\ge0.$$  Therefore,  $$(\delta'n')^2.g(x_1,\ldots,x_5)=\frac{1}{2}\left(\sum\limits_{i = 3}^5 {|I_i|} \right)\delta'n'+\left(\sum\limits_{i = 1}^5 {|I_i||I_{i+2}|} \right),$$ and $(\ref{5-2})$ implies that $(*)$.%  furthermore, the $G'[X_6]$ of the extremal graph distribute as follow: $(i)$: $I_i=\frac{\delta'n'}{2}$. $(ii)$: all edges but $G'[I_1,I_2,I_3,I_4,I_5]$ are color $1$. $(iii)$:$G_2'[I_1,I_2,I_3,I_4,I_5]$ be a complete $5$-partite $K_{\frac{\delta'n'}{2},\ldots, \frac{\delta'n'}{2}} $. $(iv)$: $G_1'[X_6]$ is $\delta'n'$-regular.
\medskip

Now we show $(**)$ in the end of Case 2  of the proof of Theorem \ref{zhu-2}.

\medskip
Let $x_i=|I_i|/(\delta'n')$ and $y_i=|Z_{i+3}|/(\delta'n')$. From Claim \ref{4-1-16},  for $i\in [5]$, $$|Z_i|+|I_{i+2}|+|I_{i+3}| \le \delta'n'.$$ Then for all $i\in [5]$, $$x_i+x_{i+1}+y_i\le 1,~ x_i\ge0 ,~  y_i\ge0.$$ Therefore,  $$(\delta'n')^2.f(x_1,\ldots,x_5,y_1,\ldots,y_5)=  \frac{1}{5}|Z|\delta'n'+ \frac{3}{10}|I|\delta'n'+\sum\limits_{i = 1}^5 {{|Z_i|}| {{I_i}} |}+\sum_{i=1}^{5}|I_i||I_{i+2}|,$$ and $(\ref{5-1})$ implies that $(**)$. %when $x_i=0.5$ and $y_i=0$ for all $i\in[5]$, $e(G')$ get the maximum value $\frac{5}{12}n'^2+\frac{1}{2}\delta'n'^2+\frac{5}{2}\delta'^2n'^2$.%However, the $G'[X_6]$ of the extremal graph distribute also as follow: $(i)$: $I_i=\frac{\delta'n'}{2}$. $(ii)$: all edges but $G'[I_1,I_2,I_3,I_4,I_5]$ are color $1$. $(iii)$:$G_2'[I_1,I_2,I_3,I_4,I_5]$ be a complete $5$-partite $K_{\frac{\delta'n'}{2},\ldots, \frac{\delta'n'}{2}} $. $(iv)$: $G_1'[X_6]$ is $\delta'n'$-regular. Thus, combine with Claim \ref{b-e-1}, it also have $e(G')\le\frac{5}{12}n'^2+\frac{1}{2}\delta'n'^2+\frac{5}{2}\delta'^2n'^2$.

%The code to calculate $\mathop {\max }\limits_{{D_f}} f$ and $\mathop {\max }\limits_{{D_g}} g$ using the LINGO are as follows

\medskip
The codes for calculating $\mathop {\max }\limits_{{D_f}} f$ are as follows:

\noindent
max=0.3*(x1+x2+x3+x4+x5)+0.2*(y1+y2+y3+y4+y5)\\+x1*(x3+x4)+x2*(x4+x5)+x3*x5+x1*y3+x2*y4+x3*y5+x4*y1+x5*y2;\\
$x1+x2+y1<1;\\
x2+x3+y2<1;\\
x3+x4+y3<1;\\
x4+x5+y4<1;\\
x5+x1+y5<1;\\
0<x1;\\
0<x2;\\
0<x3;\\
0<x4;\\
0<x5;\\
0<y1;\\
0<y2;\\
0<y3;\\
0<y4;\\
0<y5;$
\medskip

The codes  for calculating $\mathop {\max }\limits_{{D_g}} g$ are as follows:

\noindent
max=0.5*(x3+x4+x5)+x1*(x3+x4)+x2*(x4+x5)+x3*x5;\\
$x1+x2<1;\\
x2+x3<1;\\
x3+x4<1;\\
x4+x5<1;\\
x5+x1<1;\\
0<x1;\\
0<x2;\\
0<x3;\\
0<x4;\\
0<x5;$

\end{spacing}

\begin{thebibliography}{99}
\bibitem{afk} N.~Alon, E.~Fischer, M.~Krivelevich and M.~Szegedy, Efficient testing of large graphs, {\em Combinatorica} {20} (2000), 451--476.
 \bibitem{bhs}
J. Balogh, P. Hu, and M. Simonovits, Phase transitions in the Ramsey-Tur\'an theory, {\em J. Combin. Theory Ser. B} 114 (2015), 148--169.

 \bibitem{bls}
J. Balogh, H. Liu, and M. Sharifzadeh, On two problems in Ramsey-Tur\'an theory, {\em SIAM J. Discrete Math} 31 (2017), 1848--1866.

 \bibitem{B1}
  B. Bollob\'{a}s and P. Erd\H{o}s, On a Ramsey-Tur\'{a}n type problem, {\em J. Combin. Theory Ser. B} 21 (1976), 166--168.

% \bibitem{bm}
%J. Bondy and U. Murty, Graph Theory with Applications, Springer, 2008.

  \bibitem{Brandt}
  S.Brandt, Triangle-free graphs whose independence number equals the degree, {\em Discrete Math.} 310 (2010), 662--669.

 \bibitem{B2}
  W. G. Brown, P. Erd\H{o}s and M. Simonovits, Extremal problems for directed graphs, {\em J. Combin. Theory Ser. B} 15 (1973), 77--93.
 \bibitem{BLS}
J. Balogh, H. Liu and M. Sharifzadeh, On two problems in Ramsey-Tur\'{a}n theory, {\em SIAM J. Discrete math.} 31 (2017), 1848--1866.

\bibitem{cf}
M. Codish, M. Frank, A. Itzhakov and A. Miller, Computing the Ramsey Number $R(4,3,3)$ Using
Abstraction and Symmetry Breaking, {\em Constraints} 21 (2016), 375--393.

\bibitem{Conlon}
D. Conlon, The Ramsey number of books, {\em Adv. Combin.} 3 (2019), 12pp.
\bibitem{c-f}
D. Conlon, J. Fox and Y. Wigderson, Ramsey number of books and quasirandomness, {\em Combinatorica} 42 (2022), no. 3, 309--363.

  \bibitem{erd}
  P. Erd\H{o}s, Graph theory and probability II, {\em Canad. J. Math.} 13 (1961), 346--352.
\bibitem{E-H-3S}
  P. Erd\H{o}s, A. Hajnal, M. Simonovits, V. S\'{o}s and E. Szemer\'{e}di, Tur\'{a}n-Ramsey theorems and $K_p$-independnece numbers, \emph{Combin. Probab. Comput.} 3 (1994), 297--325.
  \bibitem{B3}
  P. Erd\H{o}s, A. Hajnal, M. Simonovits, V. S\'{o}s and E. Szemer\'{e}di, Tur\'{a}n-Ramsey theorems and simple asymptotically extremal structures, {\em Combinatorica} 13 (1993), 31--56.
  \bibitem{B4}
  P. Erd\H{o}s, A. Hajnal, V. S\'{o}s and E. Szemer\'{e}di, More results on Ramsey-Tur\'{a}n type problems, {\em Combinatorica} 3 (1983), 69--81.
  \bibitem{B5}
   P. Erd\H{o}s and C. A. Rogers, The construction of certain graphs, {\em Canadian J. Math} 14 (1962), 702--707.
  \bibitem{B6}
   P. Erd\H{o}s and V. T. S\'{o}s, some remarks on Ramsey's and Tur\'{a}n's theorem, in Combinatorial Theory and Its Applications, II (Proc. Colloq., Balatonf\"{u}red, 1969), North-Holland, Amsterdam, 1970, 395--404.

   \bibitem{e70}
 P. Erd\H{o}s, On the graph theorem of Tur\'{a}n, {\em Mat. Lapok} 21 (1970), 249--251 (in Hungarian).

   \bibitem{es2}
   P. Erd\H{o}s and V. T. S\'{o}s, On Tur\'{a}n-Ramsey's type theorems, II, {\em Stud. Sci. Math. Hung.} 14 (1979), 27--36.
\bibitem{Erdo-stone}
P. Erd\H{o}s and A. H. Stone, On the structure of linear graphs, \emph{ Bull. Amer. Math. Soc.} 52 (1946), 1087--1091.
  \bibitem{furedi}
  Z. F\"{u}redi, A proof of the stability of extremal graphs, Simonovits'stability from Szemer\'{e}di's regularity, {\em J. Combin. Theory Ser. B} 115 (2015), 66--71.


     \bibitem{flz}
 J. Fox, P. Loh, and Y. Zhao, The critical window for the classical Ramsey-Tur\'an problem, {\em Combinatorica} 35 (2015), 435--476.


     \bibitem{fs}
      J. Fox and B. Sudakov, Dependent random choice, {\em Random Structures Algorithms} 38 (2011), 68--99.

    \bibitem{gg}
R. E. Greenwood and A. M. Gleason, Combinatorial relations and chromatic graphs, {\em Canad. J. Math.} 7 (1955), 1--7.

 \bibitem{H-L}
X. Hu and Q. Lin, Two-colored Ramsey-Tur\'{a}n densities involving triangles, {\em SIAM J. Discrete Math.} 38 (2024), 2132--2162.

  \bibitem{kim}
J. H. Kim, The Ramsey number $R(3, t)$ has order of magnitude $t^2/ \log t$, {\em Random Structures Algorithms} 7 (1995), 173--207.

  \bibitem{kkl}
  J. Kim, Y. Kim and H. Liu, Two conjectures in Ramsey-Tur\'{a}n theory, {\em SIAM J. Discrete Math.} 33 (2019), 564--586.

\bibitem{ko-sim}
J. Koml\'{o}s and M. Simonovits, Szemer\'{e}di's regularity lemma and its applications to graph theory.
{\em Combinatorics, Paul Erd\H{o}s is eighty, Vol. 2 (Keszthely, 1993)}, 295--352, Bolyai Soc. Math. Stud., 2, {\em J\'{a}nos Bolyai Math. Soc., Budapest}, 1996.

%J. Koml\'{o}s and M. Simonovits, Szemer\'{e}di's regularity lemma and its applications to graph theory, {\em Combinatorica}, 1996, 295--352.
\bibitem{li-l} Y. Li and Q. Lin, Elementary methods of graph Ramsey theory, Springer, 2022.


\bibitem{ll}
M. Liu and Y. Li, Two results on Ramsey-Tur\'{a}n numbers, {\em Electron. J. Combin.} 28 (2021), no. 4, Paper No. 4.6, 8 pp.

  \bibitem{lr}
C. M. L\"{u}ders and C. Reiher, The Ramsey-Turan problem for cliques, {\em Israel J. Math.} 230 (2019), no. 2, 613--652.

  \bibitem{rad}
S. Radziszowski, Small Ramsey numbers, Electron. J. Combin. 1 (1994), a dynamic survey.


 \bibitem{ram}
F.~P.~Ramsey, On a problem of formal logic,  {\em Proc.~Lond. Math.~Soc.} 30 (1929), 264--286.
   \bibitem{ss}
 M. Simonovits and V. S\'os, Ramsey-Tur\'an theory, {\em Discrete Math.} 229 (2001), 293--340.
    \bibitem{sud}
 B. Sudakov, A few remarks on the Ramsey-Tur\'an-type problems, {\em J. Combin. Theory Ser. B} 88 (2003), 99--106.
\bibitem{sze78}
E.~Szemer\'edi, {Regular partitions of graphs}, in: Probl\`emes Combinatories et th\'eorie des graphs, Colloque Inter. CNRS, Univ. Orsay, Orsay, 1976, J. Bermond, J. Fournier,  M. Las Vergnas, and D. Scotteau, Eds. (1978), pp. 399--402.

  \bibitem{B8}
   E. Szemer\'{e}di, On graphs containing no complete subgraph with $4$ vertices, {\em Mat. Lapok} 23 (1972), 113--116.
 \bibitem{B9}
 P. Tur\'{a}n, Eine Extremalaufgabe aus der Fraphentheorie, {\em Mat. Fiz. Lapok} 48 (1941), 436--452.
  \bibitem{turan54}
 P. Tur\'{a}n, On the theory of graphs, {\em Colloq. Math.} 3 (1954), 19--30.
\end{thebibliography}
\end{document}